\documentclass{article}
\usepackage{amsmath,amsthm}
\usepackage{verbatim,amsmath,amscd,amssymb,amsthm}
\begin{document}
\font\germ=eufm10
\def\ssl{\hbox{\germ sl}}
\def\slh{\widehat{\ssl_2}}
\def\ge{\hbox{\germ g}}
\def\aaa{@}
\def\aaa{@}
\title{\Large\bf Nilpotent Representations of \\
Classical Quantum Groups at Roots of Unity }

\author{
Yuuki A\textsc{be} 
\thanks
{
e-mail: yu-abe@hoffman.cc.sophia.ac.jp
}
\\
Department of Mathematics, 
\\
Sophia University
\and
Toshiki N\textsc{akashima} 
\thanks
{
supported in part by JSPS Grants in Aid for 
Scientific Research, e-mail: toshiki@mm.sophia.ac.jp
}
\\
Department of Mathematics, 
\\
Sophia University,
}
\date{}

\maketitle
\centerline{\it In memory of W.A.Schnizer}

\begin{abstract}
Properly specializing the parameters in ``Schnizer modules'', 
for type A,B,C and D,
we get its unique primitive vector. Then we show that
the module generated by the primitive vector is an
irreducible highest weight module of 
finite dimensional classical 
quantum groups at roots of unity.
\end{abstract}
\maketitle

\renewcommand{\labelenumi}{$($\roman{enumi}$)$}
\renewcommand{\labelenumii}{$(${\rm \alph{enumii}}$)$}
\font\germ=eufm10

\def\al{\alpha}
\def\d{\delta}
\def\deg{\textrm{deg}}
\def\e{\varepsilon}
\def\g{\texttt{g}}
\def\lm{\lambda}
\def\nn{\nonumber}
\def\ot{\otimes}
\def\ue{U_{\varepsilon}}
\def\uf{U_{\varepsilon} ^{\rm fin}}
\def\uq{U_q(\texttt{g})}
\def\ur{U_{\varepsilon} ^{\rm res}}
\def\q{\quad}
\def\qq{\qquad}
\def\vp{\varphi}
\def\vep{\varepsilon}
\def\w{{\mathcal W}}
\def\z{\bbZ}

\newcommand{\bm}{{\bf m}}
\newcommand{\bU}{{\bf U}}
\newcommand{\cI}{{\mathcal I}}
\newcommand{\cA}{{\mathcal A}}
\newcommand{\cB}{{\mathcal B}}
\newcommand{\cC}{{\mathcal C}}
\newcommand{\cD}{{\mathcal D}}
\newcommand{\cF}{{\mathcal F}}
\newcommand{\cH}{{\mathcal H}}
\newcommand{\cK}{{\mathcal K}}
\newcommand{\cL}{{\mathcal L}}
\newcommand{\cM}{{\mathcal M}}
\newcommand{\cMod}{{\mathcal M}\!{\it od}}
\newcommand{\cN}{{\mathcal N}}
\newcommand{\cO}{{\mathcal O}}
\newcommand{\cS}{{\mathcal S}}
\newcommand{\cV}{{\mathcal V}}
\newcommand{\fra}{\mathfrak a}
\newcommand{\frb}{\mathfrak b}
\newcommand{\frc}{\mathfrak c}
\newcommand{\frd}{\mathfrak d}
\newcommand{\fre}{\mathfrak e}
\newcommand{\frf}{\mathfrak f}
\newcommand{\frg}{\mathfrak g}
\newcommand{\frh}{\mathfrak h}
\newcommand{\fri}{\mathfrak i}
\newcommand{\frj}{\mathfrak j}
\newcommand{\frk}{\mathfrak k}
\newcommand{\frI}{\mathfrak I}
\newcommand{\fm}{\mathfrak m}
\newcommand{\frn}{\mathfrak n}
\newcommand{\frp}{\mathfrak p}
\newcommand{\fq}{\mathfrak q}
\newcommand{\frr}{\mathfrak r}
\newcommand{\frs}{\mathfrak s}
\newcommand{\frt}{\mathfrak t}
\newcommand{\fru}{\mathfrak u}
\newcommand{\frA}{\mathfrak A}
\newcommand{\frB}{\mathfrak B}
\newcommand{\frF}{\mathfrak F}
\newcommand{\frG}{\mathfrak G}
\newcommand{\frH}{\mathfrak H}
\newcommand{\frJ}{\mathfrak J}
\newcommand{\frN}{\mathfrak N}
\newcommand{\frP}{\mathfrak P}
\newcommand{\frT}{\mathfrak T}
\newcommand{\frU}{\mathfrak U}
\newcommand{\frV}{\mathfrak V}
\newcommand{\frX}{\mathfrak X}
\newcommand{\frY}{\mathfrak Y}
\newcommand{\frZ}{\mathfrak Z}
\newcommand{\rA}{\mathrm{A}}
\newcommand{\rC}{\mathrm{C}}
\newcommand{\rd}{\mathrm{d}}
\newcommand{\rB}{\mathrm{B}}
\newcommand{\rD}{\mathrm{D}}
\newcommand{\rE}{\mathrm{E}}
\newcommand{\rH}{\mathrm{H}}
\newcommand{\rK}{\mathrm{K}}
\newcommand{\rL}{\mathrm{L}}
\newcommand{\rM}{\mathrm{M}}
\newcommand{\rN}{\mathrm{N}}
\newcommand{\rR}{\mathrm{R}}
\newcommand{\rT}{\mathrm{T}}
\newcommand{\rZ}{\mathrm{Z}}
\newcommand{\bbA}{\mathbb A}
\newcommand{\bbC}{\mathbb C}
\newcommand{\bbG}{\mathbb G}
\newcommand{\bbF}{\mathbb F}
\newcommand{\bbH}{\mathbb H}
\newcommand{\bbP}{\mathbb P}
\newcommand{\bbN}{\mathbb N}
\newcommand{\bbQ}{\mathbb Q}
\newcommand{\bbR}{\mathbb R}
\newcommand{\bbV}{\mathbb V}
\newcommand{\bbZ}{\mathbb Z}
\newcommand{\adj}{\operatorname{adj}}
\newcommand{\Ad}{\mathrm{Ad}}
\newcommand{\Ann}{\mathrm{Ann}}
\newcommand{\rcris}{\mathrm{cris}}
\newcommand{\ch}{\mathrm{ch}}
\newcommand{\coker}{\mathrm{coker}}
\newcommand{\diag}{\mathrm{diag}}
\newcommand{\Diff}{\mathrm{Diff}}
\newcommand{\Dist}{\mathrm{Dist}}
\newcommand{\rDR}{\mathrm{DR}}
\newcommand{\ev}{\mathrm{ev}}
\newcommand{\Ext}{\mathrm{Ext}}
\newcommand{\cExt}{\mathcal{E}xt}
\newcommand{\fin}{\mathrm{fin}}
\newcommand{\Frac}{\mathrm{Frac}}
\newcommand{\GL}{\mathrm{GL}}
\newcommand{\Hom}{\mathrm{Hom}}
\newcommand{\hd}{\mathrm{hd}}
\newcommand{\rht}{\mathrm{ht}}
\newcommand{\id}{\mathrm{id}}
\newcommand{\im}{\mathrm{im}}
\newcommand{\inc}{\mathrm{inc}}
\newcommand{\ind}{\mathrm{ind}}
\newcommand{\coind}{\mathrm{coind}}
\newcommand{\Lie}{\mathrm{Lie}}
\newcommand{\Max}{\mathrm{Max}}
\newcommand{\mult}{\mathrm{mult}}
\newcommand{\op}{\mathrm{op}}
\newcommand{\ord}{\mathrm{ord}}
\newcommand{\pt}{\mathrm{pt}}
\newcommand{\qt}{\mathrm{qt}}
\newcommand{\rad}{\mathrm{rad}}
\newcommand{\res}{\mathrm{res}}
\newcommand{\rgt}{\mathrm{rgt}}
\newcommand{\rk}{\mathrm{rk}}
\newcommand{\SL}{\mathrm{SL}}
\newcommand{\soc}{\mathrm{soc}}
\newcommand{\Spec}{\mathrm{Spec}}
\newcommand{\St}{\mathrm{St}}
\newcommand{\supp}{\mathrm{supp}}
\newcommand{\Tor}{\mathrm{Tor}}
\newcommand{\Tr}{\mathrm{Tr}}
\newcommand{\wt}{\mathrm{wt}}
\newcommand{\Ab}{\mathbf{Ab}}
\newcommand{\Alg}{\mathbf{Alg}}
\newcommand{\Grp}{\mathbf{Grp}}
\newcommand{\Mod}{\mathbf{Mod}}
\newcommand{\Sch}{\mathbf{Sch}}\newcommand{\bfmod}{{\bf mod}}
\newcommand{\Qc}{\mathbf{Qc}}
\newcommand{\Rng}{\mathbf{Rng}}
\newcommand{\Top}{\mathbf{Top}}
\newcommand{\Var}{\mathbf{Var}}
\newcommand{\gromega}{\langle\omega\rangle}
\newcommand{\lbr}{\begin{bmatrix}}
\newcommand{\rbr}{\end{bmatrix}}

\newcommand{\forb}{\bigcirc\kern-2.8ex \because}
\newcommand{\forbb}{\bigcirc\kern-3.0ex \because}
\newcommand{\forbbb}{\bigcirc\kern-3.1ex \because}
\newcommand{\SpS}{spectral sequence}
\newcommand\C{\mathbb C}
\newcommand\hh{{\hat{H}}}
\newcommand\eh{{\hat{E}}}
\newcommand\F{\mathbb F}
\newcommand\fh{{\hat{F}}}
\newcommand{\End}{\operatorname{End}}
\newcommand{\Stab}{\operatorname{Stab}}
\newcommand{\mo}{\operatorname{mod}}
\newcommand\pf{\noindent {\bf Proof:  }}

\def\al{\alpha}
\def\c{\textbf{C}}
\def\cq{\textbf{C}(q)}
\def\d{\delta}
\def\deg{\textrm{deg}}
\def\e{\varepsilon}
\def\g{\frg}
\def\l{\lambda}
\def\no{\nonumber}

\section{Introduction}
\renewcommand{\thesection}{\arabic{section}}
\setcounter{equation}{0}
\renewcommand{\theequation}{\thesection.\arabic{equation}}

\theoremstyle{definition}
\newtheorem{df}{Definition}[section]
\theoremstyle{plain}
\newtheorem{pro}[df]{Proposition}
\newtheorem{lem}[df]{Lemma}
\newtheorem{thm}[df]{Theorem}
\newtheorem{ex}[df]{Example}

The representation theory of quantum groups 
at roots of unity are divided into the 
following two types: one is for $\ue$ defined by DeConcini-Kac
(=non-restricted type)
\cite{DK}
and the other is for 
$\ur$ defined by Lusztig(=restricted type) \cite{L}. 
In the latter case, the classification of irreducible modules
is same as generic case, that is, they are classified by 
highest weights(\cite{L},\cite{L2}). 
In the former case, however, 
most irreducible modules are no longer highest or
lowest weight modules and they are characterized by several 
continuous parameters(\cite{DK}). For type A, such modules 
are constructed very explicitly in \cite{DJMM}, 
which is called maximal cyclic representations.
For any simple Lie algebra,
Schnizer introduced an
alternative construction of such modules
in \cite{S},\cite{S2}, which we also call a maximal cyclic representation
or  ``Schnizer module''.

In \cite{N}, the second author found that for type $A_n$-case if 
the continuous parameters in maximal cyclic representations 
are specialized 
properly, then there exists a  unique primitive vector and 
the submodule generated by the primitive vector is irreducible
as a module of finite dimensional quantum group 
at roots of unity(denoted by $\uf$).
In this paper, we shall show that this method is applicable 
to the Schnizer modules
of types $A_n$, $B_n$, $C_n$ and $D_n$.

In order to explain what we shall do in this article, let us see 
$A_n$-case explicitly:
Let $N=\frac{1}{2}n(n+1)$ be the number of positive roots, $l$ be an odd integer 
greater than 3 and $\vep$
be the primitive $l$-th root of unity. Set $V:=(\bbC^l)^{\ot N}$
and for each $a,b\in(\bbC^\times)^N$ and $\lm\in \bbC^n$,
we can define a $U_\vep(\texttt{sl} (n+1, \bbC))$-module structure on $V$ 
as follows. (Indeed, the module as below is similar to the 
maximal cyclic representation as in \cite{DJMM}.).
\begin{thm}[Schnizer module \cite{S2}]
\label{sl-rep} 
For any 
$a=(a_{i,j})_{1\leq i\leq j\leq n}
\in(\bbC^{\times})^{N}$,
$b=(b_{i,j})_{1 \leq i \leq j \leq n} \in \bbC^{N}$,
$\l=(\l_1, \cdots, \l_n) \in \bbC^n$, 
we obtain a $U_{\e}(\texttt{sl} (n+1, \bbC))$-module structure on $V$:
$\Phi_{\l,a,b}: U_{\e}(\texttt{sl} (n+1, \bbC)) \longrightarrow 
\textrm{End}(V)$. For any $i \in I$, 
\begin{eqnarray*}
&& \Phi _{\l,a,b} (t_i) (u(m))=\e^{\mu_{i,n}^m}u(m), \\
&&\Phi _{\l,a,b} (f_i)(u(m)) 
=\sum _{k=i}^{n} [m_{i,k}-m_{i+1,k}+b_{i,k}-b_{i+1,k}-\mu_{i,k-1}^m+1]
a_{i,k}^{-1}u(m+\e_{i,k}), 
\end{eqnarray*}
\begin{eqnarray*}
&&\Phi _{\l,a,b} (e_i) (u(m))
=\sum _{k=1}^{i} [m_{k-1,n-i+k}-m_{k,n-i+k}+b_{k-1,n-i+k}-b_{k,n-i+k}+1]\\
&&\qq \qq \qq \qq \times(\prod _{p=k+1}^{i}a_{p-1,n-i+p}^{-1}a_{p,n-i+p})
u(m+\sum_{p=k+1}^i(\e_{p-1,n-i+p}-\e_{p,n-i+p})), 
\end{eqnarray*}
where for any $(i-1 \leq j \leq n)$, set 
\begin{eqnarray*}
\mu_{i,j}^m=\l_i+m_{i-1,i-1}+b_{i-1,i-1}
+\sum_{p=i}^j(m_{i-1,p}-2m_{i,p}+m_{i+1,p}+b_{i-1,p}-2b_{i,p}+b_{i+1,p}),
\end{eqnarray*}
and $a_{i,j}:=1, b_{i,j}:=0, m_{i,j}:=0$ 
if the index $(i,j)$ is out of range.
If $j>i$ then $\sum_{k=j}^i(\cdots):=0, \prod_{p=j}^i(\cdots):=1$.
\end{thm}


Here if we specialize $a,b,\lm$ as follows, we can find the 
unique primitive vector $u_\lm$ in $V$ and it has a weight
$\lm$:
Let us define  $a^{(0)}=(a^{(0)}_{i,j})_{1 \leq i \leq j \leq n} 
\in (\bbC^{\times})^{\frac{1}{2}n(n+1)},
b^{(0)}=(b^{(0)}_{i,j})_{1 \leq i \leq j \leq n} 
\in \bbC^{\frac{1}{2}n(n+1)}$ by 
\begin{eqnarray*}
a^{(0)}_{i,j}:=1, \q b^{(0)}_{i,j}:=i \q (1 \leq i \leq j \leq n).
\end{eqnarray*}
\begin{pro}
\label{sl-pv}
For any $\l = (\l _1 , \cdots , \l_n) \in \bbC ^n$, 
let $(\Phi_{\l,a^{(0)},b^{(0)}}, V)$ 
be the representation as in Theorem \ref{sl-rep}.
A vector $u\in V$ 
satisfies that 
$\Phi_{\l,a^{(0)},b^{(0)}}(e _i) u=0$ for any $i \in I$
if and only if $u \in \bbC u(0)$.
\end{pro}
Finally, it turns out that 
the submodule $U_\vep u_\lm\subset V$ is an
irreducible highest weight $\uf$-module. By this method, 
we obtain all finite dimensional irreducible $\uf$-modules:
\begin{thm}
\label{sl-nil}
For any $\l =(\l_1 , \cdots , \l_n) \in \z^n_l$
($\z_l:=\{0,1, \cdots , l-1\}$),
we define $\l^{'}:=(\l_1^{'}, \cdots , \l_n^{'}) \in \z^n$ by
\[ \l _i^{'}:=\l_i+2 \q (1 \leq i \leq n).\]     
Let $(\Phi _{\l^{'},a^{(0)},b^{(0)}}, V)$ be the 
$\ue(\texttt{sl}(n+1,\bbC))$-representation as in 
Theorem \ref{sl-rep}. 
Let $\ue u(0)$  be the $\ue(\texttt{sl}(n+1,\bbC))$-submodule of $V$ 
generated by $u(0)$.
Then $\ue u(0)$ is a finite dimensional irreducible 
$\uf(\texttt{sl}(n+1,\bbC))$-module of type $1$ with highest weight $\l$. 
\end{thm}
\noindent 
The proofs of the above statements are done by the
similar way to the ones in \cite{N}.

The organization of the paper is as follows:
in Sect 2, we prepare notations and review the theory of 
quantum groups
at roots of unity briefly. 
In Sect.3, we introduce Schnizer modules
and show the uniqueness of primitive vectors in it under 
some specialization of the parameters. In the last section, 
we show that the submodule generated by the primitive vector
is regarded as a module for the 
finite dimensional quantum group 
at roots of unity $\uf$ of types (A), B, C and D. At last,
we obtain 
 that such submodule is an irreducible $\uf$-module and 
all finite dimensional irreducible 
$\uf$-modules are exhausted by 
such modules.

\section{Quantum enveloping algebra $U_q(\frg)$ }
\renewcommand{\thesection}{\arabic{section}}
\setcounter{equation}{0}
\renewcommand{\theequation}{\thesection.\arabic{equation}}

\subsection{Definition of quantum enveloping algebra}
\q In this subsection, we define the quantum enveloping  algebra $\uq$ 
for a generic $q$.\\
\q Let $\bbC(q)$ be the rational function field in an indeterminate
 $q$. Define  
\begin{eqnarray}
& [a]_{q^d}:=\displaystyle \frac{q^{da}-q^{-da}}{q^d-q^{-d}},
\q [a]:=[a]_{q}, \nonumber \\ 
& [a]_{q^d}!:=[a]_{q^d} [a-1]_{q^d} \cdots [1]_{q^d}, 
\q [0]!:=1, \nonumber
\end{eqnarray}
for any $a,d \in \bbZ_+ :=\{0,1,2, \cdots \}.$
Let  $\frg$ be a finite dimensional simple Lie algebra 
over $\bbC$ of rank $n$ and  
 $ \{\alpha_1,\cdots ,\alpha_n \}$ be the set of simple roots, 
$I:=\{1,2, \cdots, n\}$,
 $\Delta$ be the set of roots (resp.   
$\Delta_+$ be the set of positive roots).
Define the root lattice $Q=\bigoplus_{i=1}^n \bbZ \alpha_i$
(resp.$Q_+=\bigoplus_{i=1}^n \bbZ_+ \alpha_i$). 
 Let $(\texttt{a}_{ij})_{i,j=1}^n$ be the Cartan matrix associated with $\g$,
and $d=(d_1,\cdots ,d_n)$ be an element in $\bbN^n$ 
such that $d_i \texttt{a}_{ij}=d_j\texttt{a}_{ji}$ for any $ i,j \in I$  and 
g.c.d $(d_1, \cdots ,d_n) =1$.
We denote the Weyl group of $\g$ by  $W$ which is generated by the
 simple reflections  $\{s_1, \cdots, s_n\}$. 
Now, we define the quantum enveloping algebra $\uq$
 over $\bbC (q)$.
\begin{df}
\label{1.1}
Quantum enveloping algebra $U_q(\frg)$ is an associative 
$\bbC (q)$-algebra 
generated by $\{e_i , f_i , t_i^{\pm1}|i \in I\}$
with the  relations 
\begin{eqnarray}
& t_i t_i^{-1} =t_i^{-1} t_i=1, \quad t_i t_j = t_j t_i,\nonumber \\
& t_i e_j t_i^{-1}=q_i^{\texttt{a}_{ij}}e_j, \nonumber \\
& t_i f_j t_i^{-1}=q_i^{-\texttt{a}_{ij}}e_j,\nonumber \\ 
& e_i f_j- f_j e_i = \delta_{ij} \{t_i\}_{q_i},\nonumber \\
& \sum_{k=0}^{1-\texttt{a}_{ij}} (-1)^k e_i^{(k)} 
e_j e_i^{(1-\texttt{a}_{ij}-k)}=
 \sum_{k=0}^{1-\texttt{a}_{ij}} (-1)^k f_i^{(k)} f_j 
f_i^{(1-\texttt{a}_{ij}-k)}=0
\q (i \neq j),\nonumber 
\end{eqnarray}
where 
$q_i:=q^{d_i}$, 
$e_i^{(k)}:= \displaystyle \frac{1}{[k]_{q^{d_i}}!} e_i^k,
f_i^{(k)} := \displaystyle \frac{1}{[k]_{q^{d_i}}!} f_i^k $,  
  $\{t_i\}_{q_i} := \displaystyle \frac{t_i-t_i^{-1}}{q_i-q_i^{-1}}$. \\
\q Let $U_q^+(\frg)$ (resp. $U_q^-(\frg),U_q^0(\frg)$)
 be the $\bbC(q)$-subalgebra of $\uq$ generated by
$\{e_i\}_{i=1}^n$ (resp.
$\{f_i\}_{i=1}^n,\{t_i^{\pm 1}\}_{i=1}^n$).
\end{df}

\subsection{Non-restricted specialization}
 In this subsection, we define the non-restricted specializations
$\ue$ for a root of unity $\e$. 
\begin{df}
Let $A:= \bbC [q,q^{-1}]$ be the Laurent polynomial ring,
$U_A$ be the $A$-subalgebra of $\uq$ generated by 
$\{e_i,f_i,t_i^{\pm 1},\{t_i\}_{q_i}\}_{i=1}^n$,
$l$ be an odd integer greater than $3$, and
$\e$ be a primitive $l$-th root of unity  
 such that $\e^{2d_i} \neq 1$  for any $i \in I$.
We regard $\bbC$ as $A$-algebra by 
$f(q)c:=f(\e) \cdot c$ for any $f(q) \in A, c \in \bbC$ 
and we denote it by $\bbC _{\e}$. Now we define  
\[
 \ue:=U_A \otimes _A \bbC _{\e},  
\]
and we call $\ue$ ``non-restricted specialization of $\uq$''.
By the similar manner to Definition \ref{1.1},
we define $U_{\e}^+, U_{\e}^-$ and $U_{\e}^0$, 
and we denote $u \otimes 1$ as $u$ for any $u \in U_A$. 
\end{df}
\textbf{Remark.}(\cite{DK})
One can also describe $\ue$ in term of generators and relations. 
That is, $\ue$ is an associative $\bbC$-algebra generated by
$\{e_i , f_i , t_i^{\pm 1}\}_{i=1}^n$ with the 
relations of Definition 2.1 replacing $q$ by $\e$. 

\subsection{Root vectors}
 In this subsection, 
we introduce the root vectors and its properties.
\begin{pro}[(\cite{DK}, \cite{J}]
\label{Proposition 1.3A}
\begin{enumerate}
\item
 For any $i \in I$, there exist $\ue$-automorphism $T_i$  such that  
\begin{eqnarray}
& T_i(e_i)=-f_it_i, 
\q  T_i(e_j)= \sum_{s=0}^{-\texttt{a}_{ij}} (-1)^{s-\texttt{a}_{ij}} q_i^{-s}
e_i^{(-\texttt{a}_{ij}-s)} e_j e_i^{(s)}  \q (i \neq j ), \nonumber \\
& T_i(f_i)=-t_i^{-1} e_i ,
\q  T_i(f_j)= \sum_{s=0}^{-\texttt{a}_{ij}} (-1)^{s-\texttt{a}_{ij}}
q_i^s f_i^{(s)} f_j f_i^{(-\texttt{a}_{ij}-s)} \q  (i \neq j),
\nonumber \\
& T_i(t_j)= t_j t_i^{-\texttt{a}_{ij}}. \nonumber 
\end{eqnarray}
\item
For $w \in W$, let $w=s_{i_1} \cdots s_{i_r}$ be a 
reduced expression of $w$, and set
$T_w := T_{i_1} \cdots T_{i_r}$. 
Then $T_w$ is well-defined 
 (that is, $T_w$ does not depend on a choice of reduced expression of $w$). 
\end{enumerate}
\end{pro}
\begin{df}
\label{Definition 1.3A}
Let $w_0$ be a longest element of $W$,
$w_0= s_{i_1} \cdots s_{i_N}$ be a reduced expression of $w_0$, 
and we set 
\[ \beta_1:=\al _{i_1},  \beta_2 :=s_{i_1}(\al_{i_2}),  
\cdots ,  \beta_N:=s_{i_1}\cdots s_{i_{N-1}} (\al_{i_N}), \] 
(by the theory of the classical Lie algebra,  
  $\Delta_+=\{\beta_1, \cdots , \beta_N\}$) and
\[ e_{\beta_k}:=T_{i_1} \cdots T_{i_{k-1}} (e_{i_k}),
f_{\beta_k}:=T_{i_1} \cdots T_{i_{k-1}} (f_{i_k}) 
\q  (1 \leq k \leq N). 
\]
We call these $e_{\beta_k},f_{\beta_k}$
 ``root vectors of $\ue $ ''.
\end{df}
\begin{df}
\label{Definition 1.3B}
Set $\deg ( e_i):=\al_i, 
\deg ( f_i):=-\al_i ,\deg ( t_i):=0$. 
\end{df}
These are compatible 
with the relations of $\ue$.
Therefore, we can regard $\ue$ as 
$Q$-graded algebra and we have
\[ \ue=\bigoplus_{\al \in Q}(\ue)_{\al},\qq
(\ue)_{\al}(\ue)_{\al ^{'}} 
\subset (\ue)_{\al + \al ^{'}},\] 
for any $\al, \al ^{'} \in Q$,
where$(\ue)_{\al}:=\{u \in \ue | \deg (u)= \al \}$. 
We also use the following propositions later. 
\begin{pro}[\cite{J}]
\label{Proposition 1.3B} 
We have $e_{\al} \in U_{\e}^+ \cap (\ue)_{\al}$,
$f_{\al} \in U_{\e}^- \cap (\ue)_{- \al}$ 
\q ($\al \in \Delta_+)$. 
\end{pro}
\begin{pro}[\cite{DK}]
\label{pro1.2}
Let $Z(\ue)$ be the center of $\ue$. 
We have $e_{\al}^l , f_{\al}^l , t_i^l \in Z(\ue)$ 
for any $\al \in \Delta _+, 1 \leq i \leq n$.
\end{pro}
Next, we introduce the PBW theorem 
and the triangular decomposition.
They will be used in the subsequent sections.
Let $\{ \beta_1, \cdots, \beta_N\}$ 
be as in Definition \ref{Definition 1.3A},
then $\Delta _{+}=\{ \beta_1, \cdots, \beta_N\}$. 
\begin{thm}[\cite{DK}]
\label{Theorem 1.3A}
\begin{enumerate}
\item
$\{e_{\beta _1}^{m_1} \cdots e_{\beta _N}^{m_N} | 
 m_1, \cdots ,m_N \in \bbZ _+ \}$ is a
$\bbC$-basis of $U_{\e}^+$.  
\item
$\{f_{\beta}^{m_1} \cdots f_{\beta_N}^{m_N}  | 
 m_1, \cdots , m_N \in \bbZ _+ \}$ is a
$\bbC$-basis of $U_{\e}^-$.  
\item
$\{k_1^{m_1} \cdots k_n^{m_n}  | 
m_1, \cdots , m_n \in \bbZ_+\}$
 is a $\bbC$-basis of $U_{\e}^0.$ 
\item
Let $\phi$ be the multiplication map
$\phi : U_{\e}^- \otimes U_{\e}^0 \otimes
 U_{\e}^+\longrightarrow \ue$  
($u_- \otimes u_0 \otimes u_+   \mapsto u_- u_0  u_+$). 
Then $\phi$ is an isomorphism of $\bbC$-vector space.
\end{enumerate}
\end{thm}

\section{Primitive vectors}
We keep the settings and notations as in Sect.2.

\subsection{Schnizer modules}
In this subsection, 
we introduce the Schnizer modules
 of $\ue (\texttt{sp} (2n, \bbC )), 
\ue (\texttt{so} (2n+1, \bbC ))$
and $\ue (\texttt{so} (2n, \bbC))$. 
These representations are defined 
through the representations of the
 ``Weyl algebra''. 
\begin{df}
\label{Definition 2.1}
Let $\frg = \texttt{sp}(2n, \bbC )$ 
or $\texttt{so} (2n+1, \bbC )$  
(resp. $\frg = \texttt{so} (2n, \bbC )$), 
 $H$ be a group generated by 
$\{x^{\pm}_{i,j},z^{\pm}_{i,j} | 1 \leq i,j \leq n \}$ 
(resp. $\{x^{\pm}_{i,j},z^{\pm}_{i,j}| 
1 \leq i \leq n-1, 1 \leq j \leq n \}$ )
with relations 
\begin{eqnarray*} 
&&x_{i,j} z_{i,j} = \e z_{i,j} x_{i,j}, 
\q  x_{i,j} z_{k,l} = z_{k,l} x_{i,j} \q (k,l) \neq (i,j), \\
&& x_{i,j} x_{k,l} = x_{k,l} x_{i,j}, 
\q  z_{i,j} z_{k,l} = z_{k,l} z_{i,j}. 
\end{eqnarray*}
We set ${\mathcal W} := \bbC [H]$ (= group ring of $H$), 
and call it 
``Weyl algebra''.  
\end{df}
We use the following notations in the sequel:
\[ \{h\}_{\e ^d}:= 
\displaystyle \frac{h-h^{-1}}{\e ^d -\e ^{-d}},
\q \{h\} := \{h\}_{\e}, \] 
for any $h \in H, d \in \z (d \neq 0)$. 
\begin{pro}
\label{Proposition 2.1}
Let $\frg =\texttt{sp} (2n, \bbC)$ or 
$\texttt{so} (2n+1, \bbC)$
(resp. $\frg = \texttt{so} (2n, \bbC)$), 
$V := \bigotimes _{i,j=1}^n V_{ij}$
(resp. $V:= \bigotimes _{1 \leq i \leq n-1,1 \leq j \leq n}
 V_{ij}$), 
where $V_{i,j}=\bbC^l$.
Set $u_k^{(ij)} := (\d _{k0},\d _{k1}, \cdots \d_{k,l-1}) \in V_{ij}
(0 \leq k \leq l-1)$, where $\d _{ij}$ is the Kronecker's
 delta. 
Let $X:V_{ij}\longrightarrow V_{ij}$ be the linear map defined
by
\begin{equation}
 Xu_k^{(ij)}=u_{k-1}^{(ij)} \q (u_{-1}:=u_{l-1}),
\label{Prop2.1X1}
\end{equation} 
$X_{ij} : V \longrightarrow V$ be the linear map given by
\begin{equation}
 X_{ij} (u_{k_{11}}^{(11)} \otimes \cdots \otimes
u_{k_{ij}}^{(ij)} \otimes \cdots) := u_{k_{11}}^{(11)} \otimes
\cdots \otimes (Xu_{k_{ij}}^{(ij)}) \otimes \cdots 
\label{Prop2.1X2}
\end{equation}
(i.e. $X_{ij}$ acts only on the $(i,j)$ component),
$Z:V_{ij} \longrightarrow V_{ij}$ be the linear map given by
\begin{equation}
Zu_k^{(ij)}=\e ^{k} u_k^{(ij)}, 
\label{Prop2.1Z1} \end{equation}
and 
$Z_{ij} : V \longrightarrow V$ be the linear map given by
\begin{equation}
Z_{ij} (u_{k_{11}}^{(11)} \otimes \cdots \otimes
u_{k_{ij}}^{(ij)} \otimes \cdots ) := u_{k_{11}}^{(11)} 
\otimes \cdots \otimes (Zu_{k_{ij}}^{(ij)}) \otimes \cdots 
\label{Prop2.1Z2} \end{equation}
(i.e. $Z_{ij}$ acts only on the $(i,j)$ component). 
Then these $\{X_{ij},Z_{ij}\}$ satisfies the relations in
Definition \ref{Definition 2.1}.
\end{pro}
Let $N=n^2$ (resp. $N=n(n-1)$) be the number of 
the positive roots of $\g$, 
and $a=(a_{ij}), b=(b_{ij}) \in (\bbC ^{\times})^N$. Let 
$\psi _{ab} : \w \longrightarrow End(V)$ be the homomorphism of
 $\bbC$-algebra given by
\begin{equation}
\psi _{ab}(x_{ij})=a_{ij} X_{ij},
\q \psi _{ab}(z_{ij})= b_{ij} Z_{ij}.
\label{Prop2.1XZ3} 
\end{equation} 
Then, $\psi _{ab}$ is a well-defined 
representation of ${\mathcal W}$.  \\
\q Now, we introduce the Schnizer modules of 
$\ue (\texttt{sp} (2n, \bbC)), 
 \ue (\texttt{so} (2n+1, \bbC))$,
and $\ue (\texttt{so} (2n, \bbC))$ following \cite{S}. 
\begin{thm}[\cite{S} Theorem 3.8]
\label{Theorem 2.1A} 
For $\g =\texttt{sp} (2n, \bbC) \q (n \geq 2)$, and
$\l = (\l _1 ,\cdots , \l _n) \in \bbC ^n$.\
We define the map $\varphi _{\l}: \ue \longrightarrow \w$ by 
\begin{eqnarray*}
&&\varphi _{\l} (e_1) =F_{1,1},  \\
&&\varphi _{\l} (e_j)=(\prod _{k=1}^{j-1}D_{k,j})F_{j,j} +
\sum _{q=1}^{j-1} (\prod _{p=0}^{q-1}D_{p,j})C_{q,j},
\q (2 \leq j \leq n), \\
&& \varphi _{\l} (t_j) = T_{1,j}^{-1} \q (1 \leq i <j \leq n),
\\
&&
 \varphi _{\l} (f_1):=E_{1,1}, 
\qq \varphi _{\l} (f_j):= E_{j,j} + \sum _{i=1}^{j-1} B_{i,j}
\q (2 \leq j \leq n),
\end{eqnarray*}
where 
\begin{eqnarray}
&&C_{i,j} = \{z_{i,j}^{-1} z_{i,j-1}\} x_{i,j}+
\{z_{j,i}^{-1} z_{j+1,i}\}x_{i,j} x_{j,i} x_{i,j-1}^{-1}, 
\q (1 \leq i < j \leq n-1),\nonumber  \\
&&C_{i,n} = \{z_{i,n}^2 z_{n,i}^{-2}\}_{\e ^2}x_{i,n-1} x_{i,n} x_{n,i}^2
+ \{z_{i,n-1} z_{n,i}^{-1}\} x_{i,n-1} x_{i,n} x_{n,i}
+ \{z_{i,n-1} z_{i,n}^{-2}\}_{\e ^2} x_{i,n},\nonumber \\
&& \qq \qq (i \leq i < j=n), \nonumber\\
&&D_{i,j} = x_{i,j-1}^{-1} x_{i,j} x_{j+1,i}^{-1} x_{j,i}, 
\q (1 \leq i<j \leq n-1),\nonumber \\
&&D_{i,n} = x_{i,n-1}^{-2} x_{n,i}^2 \q (1 \leq i < n),
\qq D_{i,j} =1 \q otherwise,\nonumber \\
&&F_{j,j}= \{z_{j,j}^{-1}\} x_{j,j} \q (1 \leq j \leq n-1),
\qq F_{n,n} =\{z_{n,n}^{-2}\}_{\e ^2} x_{n,n},  \nonumber 
\end{eqnarray}
\begin{eqnarray}
&&T_{i,j}= (\prod _{k=i}^{j-1} A_{k,j}) T_{j,j} \q  
(1 \leq i<j \leq n),  \nonumber\\
&&A_{i,j}:=z_{i,j-1}^{-1} z_{i,j}^2 z_{i,j+1}^{-1} z_{j+2,i}^{-1}
z_{j+1,i}^{2} z_{j,i}^{-1} \q (1 \leq i < j \leq n-2),  \nonumber\\
&&A_{i,n-1}:=z_{i,n-2}^{-1} z_{i,n-1}^2 z_{i,n}^{-2}z_{n,i}^2
z_{n-1,i}^{-1} \q (1 \leq i<n-1),  \nonumber\\
&&A_{i,n}:= z_{i,n-1}^{-2} z_{i,n}^4 z_{n,i}^{-2},
\q (1 \leq i<n),  \nonumber\\
&&T_{j,j}:=z_{j,j}^2 z_{j,j+1}^{-1} z_{j+2,j}^{-1} z_{j+1,j}^2 
z_{j+1,j+1}^{-1} z_{j+2,j+1}^{-1} \e ^{\l _j}, 
\q (1 \leq j \leq n-2),  \nonumber\\
&&T_{n-1,n-1}:= z_{n-1,n-1}^2 z_{n-1,n}^{-2} z_{n,n-1}^2  z_{n,n}^{-2}
 \e ^{\l _{n-1}},  
\qq T_{n,n}:=z_{n,n}^4 \e ^{\l _{n}},  \nonumber
\end{eqnarray}
\begin{eqnarray}
&&B_{i,j}:=\{z_{i,j+1}^{-1} z_{j+2,i}^{-1} z_{j+1,i}^2 z_{j,i}^{-1}
 z_{i,j} T_{i+1,j} \} x_{i,j}^{-1}+
\{z_{j,i}^{-1} z_{j+1,i} T_{i+1,j}\}x_{j+1,i}^{-1} 
 \q (1 \leq i<j \leq n-2), \nonumber\\
&&B_{i,n-1}:= \{z_{i,n}^{-2} z_{i,n-1} z_{n,i}^{2} z_{n-1,i}^{-1}
 T_{i+1,n-1}\}x_{i,n-1}^{-1}+ \{z_{n-1,i}^{-1} z_{ni} T_{i+1,n-1}\}
x_{n,i}^{-1} 
\q (1 \leq i<n-1),\nonumber \\
&&B_{i,n}:=\{z_{i,n}^2 z_{n,i}^{-2} T_{i+1,n}\}_{\e ^2} x_{i,n}^{-1}
\q (1 \leq i<n),\nonumber \\
&&E_{j,j}:=\{z_{j,j} z_{j,j+1}^{-1} z_{j+2,j}^{-1} z_{j+1,j}^2
 z_{j+1,j+1}^{-1} z_{j+2,j+1}^{-1} \e ^{\l _j}\} x_{j,j}^{-1}
 + \{z_{j+1,j} z_{j+1,j+1}^{-1} z_{j+2,j+1}^{-1} \e ^{\l _j}\}
 x_{j+1,j}^{-1}, \nonumber\\
&& \qq \qq  (1 \leq j \leq n-2), \nonumber\\
&&E_{n-1,n-1}:=\{z_{n-1,n-1} z_{n-1,n}^{-2} z_{n,n-1}^{2} z_{n,n}^{-2}
 \e ^{\l _{n-1}}\}x_{n-1,n-1}^{-1}
 +\{z_{n,n-1} z_{n,n}^{-2}
 \e ^{\l _{n-1}}\} x_{n,n-1}^{-1}, \nonumber\\
&&E_{n,n}:=\{z_{n,n}^2 \e ^{\l _n}\}_{\e ^2} x_{n,n}^{-1}.\nonumber
\end{eqnarray}
Then $\varphi _{\l}$ is a homomorphism of $\bbC$-algebra. 
In particular, a pair $(\Phi _{\l ,a,b}:= \psi _{ab} \circ \varphi _{\l}, V)$
is a representation of $\ue(\texttt{sp}(2n,\bbC))$.
\end{thm}
We call the representation in the above theorem
``Schnizer module'' or  ``maximal cyclic representation''.

\textbf{Remark.}
\begin{enumerate}
\item
The explicit form of the actions of the generators above
are slightly different from those in \cite{S}. 
Through the $\bbC$-algebra $\ue$-automorphism $\omega$ to $\ue$
such that $(\omega (e_i), \omega (f_i), \omega (t_i))
= (f_i ,e_i, t_i^{-1})$, we have that 
the action of $e_j$, (resp. $f_j$, $t_j$) 
in \cite{S} corresponds to 
the action of $f_j$(resp.  $e_j$, $t_j^{-1}$)
 as above. 
\item
We call $\ue$-representations such that 
$e_i^l \neq 0$ and $f_i^l \neq 0$ for any $i \in I$ 
(resp. $e_i^l=0$ and $f_i^l=0$ for any $i$) 
``cyclic $\ue$-representations'' 
(resp. ``nilpotent $\ue$-representations''). 
In particular, we call $l^N$-dimensional irreducible cyclic
$\ue$-representations 
``maximal cyclic $\ue$-representations'' 
($l^N$ is the dimension of the representation in 
Theorem \ref{Theorem 2.1A}). 
Because the dimension of the finite dimensional irreducible 
$\ue$-representations are less than or equal to 
$l^N$ (\cite{DK}).
The representations of Theorem \ref{Theorem 2.1A}
 are not necessarily irreducible or cyclic.
However here, we also call these representations
maximal cyclic $\ue$-representations. 
\end{enumerate}
\begin{thm}[\cite{S} Theorem 3.10]
\label{Theorem 2.1B} 
For $\g = \texttt{so} (2n+1, \bbC) \q (n \geq 3)$ and
$\l =(\l _1 , \cdots , \l _n ) \in \bbC ^n$.
We define the map 
$\vp _{\l} : \ue \longrightarrow {\mathcal W}$ by,
\begin{eqnarray*}
&& \vp _{\l}(e_1):=F_{11}, 
\qq \vp_{\l}(e_j):= (\prod _{k=1}^{j-1} D_{k,j}) F_{j,j} +
\sum _{q=1}^{j-1} (\prod _{p=0}^{q-1} D_{p,j})C_{q,j},
\q (2 \leq j \leq n),\\
&& \vp_{\l}(t_j)= T_{1,j}^{-1}, \q (1 \leq j \leq n),\\
&& \vp _{\l}(f_1):=E_{1,1}, 
\qq \vp _{\l}(f_j):=E_{j,j} + \sum _{i=1}^{j-1} B_{i,j},
\q (2 \leq j \leq n),
\end{eqnarray*}
where 
\begin{eqnarray}
&&C_{i,j}:= \{z_{i,j-1}^2 z_{i,j}^{-2} \}_{\e ^2} x_{i,j} +
\{z_{j+1,i}^2 z_{j,i}^{-2}\}_{\e^2} x_{i,j-1}^{-1} x_{i,j} x_{j,i},
\q (1 \leq i<j \leq n-1), \nonumber\\
&&C_{i,n}:= \{z_{i,n-1}^2 z_{i,n}^{-1}\} x_{i,n} +
\{z_{n,i}^{-2} z_{i,n}\} x_{i,n-1}^{-1} x_{i,n} x_{n,i},
\q (1 \leq i<j =n),\nonumber \\
&&D_{i,j} :=x_{i,j-1}^{-1} x_{i,j} x_{j+1,i}^{-1} x_{j,i},
\q (1 \leq i<j \leq n-1),\nonumber \\
&&D_{i,n}:= x_{i,n-1}^{-1} x_{n,i}, \q (1 \leq i<n), 
\qq D_{i,j}:=1, \q otherwise, \nonumber\\
&&F_{j,j}:= \{z_{j,j}^{-2}\}_{\e ^2}x_{j,j}, \q (1 \leq j \leq n-1),
\qq F_{n,n}:= \{z_{n,n}^{-1}\} x_{n,n},
\nonumber
\end{eqnarray}
\begin{eqnarray}
&&\hspace{-20pt}T_{i,j}:=(\prod_{k=i}^{j-1} A_{k,j}) T_{j,j}, 
\q (1 \leq i <j \leq n), \nonumber\\
&&\hspace{-20pt}A_{i,j}
:= z_{i,j-1}^{-2} z_{i,j}^4 z_{i,j+1}^{-2} z_{j+2,i}^{-2}
z_{j+1,i}^{4} z_{j,i}^{-2}, \q (1 \leq i<j \leq n-2), 
\nonumber\\
&&\hspace{-20pt}
A_{i,n-1}:=z_{i,n-2}^{-2} z_{i,n-1}^4 z_{i,n}^{-2} z_{n,i}^4
 z_{n-1,i}^{-2}, \q (1 \leq i<n-1),\nonumber \\
&&\hspace{-20pt}
A_{i,n}:= z_{i,n-1}^{-2} z_{i,n}^2 z_{n,i}^{-2},
\q  (1 \leq i<n), \nonumber\\
&&\hspace{-20pt}
T_{j,j}:= z_{j,j}^4 z_{j,j+1}^{-2} z_{j+2,j}^{-2} z_{j+1,j}^4
 z_{j+1,j+1}^{-2} z_{j+2,j+1}^{-2} \e ^{\l _j},  
\q (1 \leq j \leq n-2), \nonumber\\
&&\hspace{-20pt}
T_{n-1,n-1}:=
z_{n-1,n-1}^4 z_{n-1,n}^{-2} z_{n,n-1}^4 z_{n,n}^{-2}
 \e ^{\l _{n-1}},
\qq T_{n,n}:=z_{n,n}^2 \e ^{\l _n},\nonumber\\
&&B_{i,j}:=\{z_{i,j}^2 z_{i,j+1}^{-2} z_{j+2,i}^{-2} z_{j+1,i}^4
 z_{j,i}^{-2} T_{i+1,j}\}_{\e ^2} x_{i,j}^{-1}
+ \{z_{j+1,i}^2 z_{j,i}^{-2} T_{i+1,j}\}_{\e ^2}x_{j+1,i}, \nonumber\\
&&\qq \qq \qq  (1 \leq i <j \leq n-2), \nonumber\\
&&B_{i,n-1}:=\{z_{i,n-1}^2 z_{i,n}^{-2} z_{n,i}^4 z_{n-1,i}^{-2}
 T_{i+1,n-1}\}_{\e ^2} x_{i,n-1}^{-1} +
\{z_{n,i}^2 z_{n-1,i}^{-2} T_{i+1,n-1}\}_{\e ^2} x_{n,i}^{-1}, \nonumber\\
&&\qq\qq \qq (1 \leq i<n-1),\nonumber \\
&&B_{i,n}:=\{z_{i,n} z_{n,i}^{-2} T_{i+1,n}\} x_{i,n}^{-1},
\q (1 \leq i<n), \nonumber\\
&&E_{j,j}:=\{z_{j,j}^2 z_{j,j+1}^{-2} z_{j+2,j}^{-2} z_{j+1,j}^4
 z_{j+1,j+1}^{-2} z_{j+2,j+1}^{-2} \e ^{\l _j}\}_{\e ^2} x_{j,j}^{-1}
\nonumber \\
&& \qq + \{z_{j+1,j}^2 z_{j+1,j+1}^{-2} z_{j+2,j+1}^{-2}
 \e ^{\l _j}\}_{\e ^2}x_{j+1,j}^{-1},  
 \q  (1 \leq j \leq n-2),\nonumber \\
&&E_{n-1,n-1}:=\{z_{n-1,n-1}^2 z_{n-1,n}^{-2} z_{n,n-1}^4 z_{n,n}^{-2}
\e ^{\l _{n-1}}\}_{\e ^2} x_{n-1,n-1}^{-1}
+\{z_{n,n-1}^2 z_{n,n}^{-2} \e ^{\l _{n-1}}\} x_{n,n-1}^{-1}, 
\nonumber \\
&&E_{n,n}:=\{z_{n,n} \e ^{\l _n}\} x_{nn}^{-1}.\nonumber
\end{eqnarray}
Then $\varphi _{\l}$ is a homomorphism of $\bbC$-algebra.  
In particular, 
a pair $(\Phi _{\l ,a,b} := \psi _{ab} \circ \vp _{\l} , V)$
is a representation of $\ue(\texttt{so}(2n+1,\bbC))$. 
\end{thm}
\begin{thm}[\cite{S} Theorem 3.11]
\label{Theorem 2.1C} 
For $\g =\texttt{so} (2n,\bbC ) \q (n \geq 4)$ and 
$\l =(\l _1 , \cdots , \l_n) \in \bbC ^n$.
We define the map 
$\vp_{\l} : \ue \longrightarrow \w$ by
\begin{eqnarray}
&&\hspace{-40pt} \vp_{\l}(e_1):=F_{1,1}, 
\qq \vp_{\l}(e_j):=(\prod_{k=1}^{j-1}D_{k,j})F_{j,j} +
\sum_{q=1}^{j-1}(\prod_{p=0}^{q-1}D_{p,j})C_{q,j}, \q
(2 \leq j \leq n-2), \nonumber\\
&&\hspace{-40pt}
\vp_{\l}(e_{n-1}):=(\prod_{k=1}^{n/2-1}D_{2k-1,n-1})
(\prod_{k^{'}=1}^{n/2-1}D_{2k^{'},n})F_{n-1,n-1} \nonumber\\
&&\hspace{-40pt}
+\sum_{q=1}^{n/2-1} (\prod_{p=0}^{q-1}D_{2p-1,n-1})
(\prod_{p^{'}=0}^{q-1}D_{2p^{'},n})C_{2q-1,n-1}
+\sum_{q=1}^{n/2-1}(\prod_{p=0}^{q}D_{2p-1,n-1})
(\prod_{p^{'}=0}^{q-1}D_{2p^{'},n})C_{2q,n},
\q(n; even) \nonumber\\
&&\hspace{-40pt}
\vp_{\l}(e_{n}):=(\prod_{k=1}^{n/2-1}D_{2k,n-1})
(\prod_{k^{'}=1}^{n/2-1}D_{2k^{'}-1,n})F_{n-1,n} \nonumber \\
&&\hspace{-40pt}
+\sum_{q=1}^{n/2-1}(\prod_{p=0}^{q-1}D_{2p,n-1})
(\prod_{p^{'}=0}^{q}D_{2p^{'}-1,n})C_{2q,n-1}
+\sum_{q=1}^{n/2-1}(\prod_{p=0}^{q-1}D_{2p,n-1})
(\prod_{p^{'}=0}^{q-1}D_{2p^{'}-1,n})C_{2q-1,n},
\q(n; even) \nonumber
\end{eqnarray}
\begin{eqnarray}
&&\hspace{-40pt}
\vp_{\l}(e_{n-1}):=(\prod_{k=1}^{(n-1)/2}D_{2k-1,n-1})
(\prod_{k^{'}=1}^{(n-3)/2}D_{2k^{'},n})F_{n-1,n} \nonumber \\
&&\hspace{-40pt}
+\sum_{q=1}^{(n-1)/2}(\prod_{p=0}^{q-1}D_{2p-1,n-1})
(\prod_{p^{'}=0}^{q-1}D_{2p^{'},n})C_{2q-1,n-1}
+\sum_{q=1}^{(n-3)/2}(\prod_{p=0}^{q}D_{2p-1,n-1})
(\prod_{p^{'}=0}^{q-1}D_{2p^{'},n})C_{2q,n}, 
\q(n; odd) \nonumber
\end{eqnarray}
\begin{eqnarray}
&&\hspace{-40pt}
\vp_{\l}(e_{n}):=(\prod_{k=1}^{(n-3)/2}D_{2k,n-1})
(\prod_{k^{'}=1}^{(n-1)/2}D_{2k^{'}-1,n})F_{n-1,n-1} \nonumber \\
&&\hspace{-40pt}
+\sum_{q=1}^{(n-3)/2}(\prod_{p=0}^{q-1}D_{2p,n-1})
(\prod_{p^{'}=0}^{q}D_{2p^{'}-1,n})C_{2q,n-1}
+\sum_{q=1}^{(n-1)/2}(\prod_{p=0}^{q-1}D_{2p,n-1})
(\prod_{p^{'}=0}^{q-1}D_{2p^{'}-1,n})C_{2q-1,n}, 
\q(n; odd) \nonumber\\
&&\vp _{\l} (t_j)=T_{1,j}^{-1} \q (1 \leq j \leq n),\nonumber
 \\
&& \vp _{\l}(f_1):=E_{1,1}, 
\qq \vp _{\l}(f_j):=E_{j,j} + \sum_{i=1}^{j-1}B_{i,j},\nonumber
\q (2 \leq j \leq n), 
\end{eqnarray}
where 
\begin{eqnarray*}
&&C_{i,j}:= \{z_{i,j-1} z_{i,j}^{-1} \} x_{i,j} +
\{z_{j+1,i} z_{j,i}^{-1}\} x_{i,j-1}^{-1} x_{i,j} x_{j,i},
\q ( \leq i <j \leq n-2), \nonumber\\
&&C_{i,n-1}:= \{z_{i,n-2} z_{i,n-1}^{-1} \} x_{i,n-1} +
\{z_{i,n} z_{n-1,i}^{-1}\} x_{i,n-2}^{-1} x_{i,n-1} x_{n-1,i}, 
\q (1 \leq i < n-1), \nonumber\\
&&C_{i,n}:= \{z_{i,n-2} z_{i,n}^{-1} \} x_{i,n} +
\{z_{i,n-1} z_{n-1,i}^{-1}\} x_{i,n-2}^{-1} x_{i,n} x_{n-1,i},
\q (1 \leq i < n), \nonumber\\
\end{eqnarray*}
\begin{eqnarray}
&&D_{ij} := x_{i,j-1}^{-1} x_{i,j} x_{j+1,i}^{-1} x_{j,i},
\q (1 \leq i<j \leq n-2), \nonumber\\
&&D_{i,n-1} := x_{i,n-2}^{-1} x_{i,n-1} x_{i,n}^{-1} x_{n-1,i},
\q 1 \leq i< n-1, \nonumber\\
&&D_{in} := x_{i,n-2}^{-1} x_{i,n} x_{i,n-1}^{-1} x_{n-1,i},
\q (1 \leq i< n), 
\qq D_{i,j} := 1,  \q otherwise, \nonumber \\
&&F_{j,j}:=\{z_{j,j}^{-1} x_{j,j}\}, \q (1 \leq j \leq n-1),
\qq F_{n-1,n}:=\{z_{n-1,n}^{-1} x_{n-1,n}\}.  \nonumber
\end{eqnarray}
\begin{eqnarray}
&&T_{i,j}:= (\prod_{k=i}^{j-1} A_{k,j})T_{j,j}, 
\q (1 \leq i<j \leq n-1), \nonumber\\
&&T_{i,n}:= (\prod_{k=i}^{n-2} A_{k,n})T_{n-1,n}, 
\q (1 \leq i<n), \nonumber\\
&&A_{i,j}:=z_{i,j-1}^{-1}z_{i,j}^{2}z_{i,j+1}^{-1}z_{j+2,i}^{-1}
z_{j+1,i}^{2}z_{j,i}^{-1},  \q (1 \leq i<j \leq n-3),\nonumber \\
&&A_{i,n-2}:=z_{i,n-3}^{-1}z_{i,n-2}^{2}z_{i,n-1}^{-1}z_{i,n}^{-1}
z_{n-1,i}^{2}z_{n-2,i}^{-1},  \q (1 \leq i< n-2),\nonumber \\
&&A_{i,n-1}:=z_{i,n-2}^{-1}z_{i,n-1}^{2}z_{n-1,i}^{-1},
 \q (1 \leq i<n-1),\nonumber \\
&&A_{i,n}:=z_{i,n-2}^{-1}z_{i,n}^{2}z_{n-1,i}^{-1},
  \q (1 \leq i<n), \nonumber
\end{eqnarray}
\begin{eqnarray}
&&T_{j,j}:=z_{j,j}^{2}z_{j,j+1}^{-1}z_{j+2,j}^{-1}z_{j+1,j}^{2}
z_{j+1,j+1}^{-1}z_{j+2,j+1}^{-1} \e^{\l_j}, 
 \q (1 \leq j \leq n-3),\nonumber \\
&&T_{n-2,n-2}:=z_{n-2,n-2}^{2}z_{n-2,n-1}^{-1}z_{n-2,n}^{-1}z_{n-1,n-2}^{2}
z_{n-1,n-1}^{-1}z_{n-1,n}^{-1} \e^{\l_{n-2}}, \nonumber\\
&&T_{n-1,n-1}:=z_{n-1,n-1}^{2} \e^{\l_{n-1}}, 
\qq T_{n-1,n}:=z_{n-1,n}^{2} \e^{\l_{n}}. \nonumber
\end{eqnarray}
\begin{eqnarray}
&&B_{i,j}:=\{z_{i,j}z_{i,j+1}^{-1}z_{j+2,i}^{-1}z_{j+1,i}^{2}
z_{j,i}^{-1} T_{i+1,j}\}x_{i,j}^{-1}, 
+\{z_{j+1,i} z_{j,i}^{-1} T_{i+1,j}\}x_{j+1,i}^{-1}, \nonumber \\
&& \qq \qq \qq (1 \leq i<j \leq n-3),  \nonumber\\
&&B_{i,n-2}:=\{z_{i,n-2}z_{i,n-1}^{-1}z_{i,n}^{-1}z_{n-1,i}^{2}
z_{n-2,i}^{-1} T_{i+1,n-2}\}x_{i,n-2}^{-1},
 +\{z_{n-1,i} z_{n-2,i}^{-1} T_{i+1,n-2}\}x_{n-1,i}^{-1}, \nonumber \\
&&\qq \qq \qq (1 \leq i< n-2),  \nonumber\\
&&B_{i,n-1}:=\{z_{i,n-1}z_{n-1,i}^{-1} T_{i+1,n-1}\}x_{i,n-1}^{-1},
\q (1 \leq i< n-1),  \nonumber\\
&&B_{i,n}:=\{z_{i,n}z_{n-1,i}^{-1} T_{i+1,n}\}x_{i,n}^{-1},
\q (1 \leq i< n),  \nonumber\\
&&E_{j,j}:=\{z_{j,j}z_{j,j+1}^{-1}z_{j+2,j}^{-1}z_{j+1,j}^{2}
z_{j+1,j+1}^{-1} z_{j+2,j+1}^{-1} \e ^{\l _j}\}x_{j,j}^{-1}, \nonumber \\
&& \qq \qq \qq+\{z_{j+1,j} z_{j+1,j+1}^{-1} z_{j+2,j+1}^{-1}
 \e ^{\l _j}\}x_{j+1,j}^{-1},
\q (1 \leq j \leq n-3), \nonumber \\
&&E_{n-2,n-2}:=\{z_{n-2,n-2}z_{n-2,n-1}^{-1}z_{n-2,n}^{-1}z_{n-1,n-2}^{2}
z_{n-1,n-1}^{-1} z_{n-1,n}^{-1} \e ^{\l _{n-2}}\}x_{n-2,n-2}^{-1} \nonumber\\
&& \qq \qq \qq+\{z_{n-1,n-2} z_{n-1,n-1}^{-1} z_{n-1,n}^{-1}
 \e ^{\l _{n-2}}\}x_{n-1,n-2}^{-1}, \nonumber\\
&&E_{n-1,n-1}:=\{z_{n-1,n-1} \e ^{\l _{n-1}}\}x_{n-1,n-1}^{-1},
\qq E_{n-1,n}:=\{z_{n-1,n} \e ^{\l _{n}}\}x_{n-1,n}^{-1}. \nonumber
\end{eqnarray}
\q Then $\varphi _{\l}$ is a homomorphism of $\bbC$-algebra.  
In particular, a pair $(\Phi _{\l ,a,b} := \psi _{ab} \circ \vp _{\l} , V)$
is a representation of $\ue(\texttt{so}(2n,\bbC))$. 
\end{thm}

\subsection{Existence and uniqueness of primitive vector in $V$}
Specializing the parameters $(a,b)$ properly,  
we show the existence and uniqueness of 
primitive vector in the Schnizer modules. \\
\q First, we fix the following notations to write down 
the action of generators of the 
$\ue$ on $(\Phi _{\l ,a,b},V)$. 
 Let $N$ be the number of positive roots. We set
\[M:=\{m=(m_{ij})_{1 \leq i \leq (N/n) , 1\leq j \leq n}
 \in \bbZ^N | 0 \leq m_{ij} \leq l-1 
\, \textrm{for any}\, i,j \}.\]
 For any $m=(m_{ij}) \in M$, 
we set $u(m):=u_{m_{11}}^{(11)} \otimes u_{m_{12}}^{(12)} \otimes
\cdots \otimes u_{m_{N/n,n}}^{(N/n,n)} \in V$, and 
$\e _{ij}:=(\d _{i1} \d _{j1}, \d _{i1} \d _{j2} , \cdots ,
\d _{i,N/n} \d _{jn}) \in M$,
where $u_k^{(ij)}$ is of Proposition 3.2, and 
$\d _{ij}$ is Kronecker's delta.
Obviously, $\{u(m) | m \in M\}$ is a $\bbC$-basis of $V$. \\
\q Next, we show that $(\Phi _{\l ,a,b},V)$ has the vectors 
which is called ``primitive vectors'' by specializing the 
parameters $(a,b)$ properly. First, we write the explicit action of
$e_j$ on $V$. 
Let us start from the $\texttt{sp}(2n, \bbC)$ case.

\subsubsection{$\texttt{sp}(2n, \bbC)$-case}
\begin{lem}
\label{Lemma 2.2A}
For $\g :=\texttt{sp} (2n,\bbC)  \q (n \geq 2),
\l :=(\l _1 , \cdots , \l _n) \in \bbC ^n$,
and any $ i,j \in I$, set 
\[ a_{ij}:=0, 
\q b_{ij}:=1-i+j \q (i \leq j),
\q b_{ij}:=2n+2-i-j \q (i>j),
\]
and $a^{(0)}=(\e^{a_{ij}})_{i,j=1}^n,
b^{(0)}=(\e^{b_{ij}})_{i,j=1}^n \in (\bbC ^{\times})^{n^2}$.
Let $(\Phi_{\l,a^{(0)},b^{(0)}}, V)$ 
be the representation as in Theorem \ref{Theorem 2.1A}.
For $u= \sum _{m \in M} c_m u(m)\in V \q (c_m \in \bbC)$,
we have 
\begin{eqnarray}
&&e_1 . u= \sum _{m \in M}c_m [-m_{11}]u(m-\e _{11}), \nonumber\\
&&e_j . u= \sum _{m \in M}c_m [-m_{jj}]u(m+ \al _j) 
+ \sum _{m \in M} \sum_{q=1}^{j-1}c_m [m_{q,j-1}-m_{qj}]
u(m+ \beta _{qj})\nonumber \\
&&\qq \qq + \sum _{m \in M} \sum_{q=1}^{j-1}c_m [m_{j+1,q}-m_{jq}]
u(m+ \beta _{qj}^{'}),  
\q   (2 \leq j \leq n-1), \nonumber
\end{eqnarray}
\begin{eqnarray}
&&e_n . u= \sum _{m \in M}c_m [-2m_{nn}]_{\e ^2}u(m+ \al _n)  
 + \sum _{m \in M} \sum_{q=1}^{n-1}c_m [2(m_{qn}-m_{nq})]_{\e ^2}
u(m+ \beta _{qn})\nonumber \\
&&\qq \q  + \sum _{m \in M} \sum_{q=1}^{n-1}c_m [m_{q,n-1}-m_{nq}]
u(m+ \beta _{qn}^{'})  
 + \sum _{m \in M} \sum_{q=1}^{n-1}c_m [2(m_{qn}-m_{nq})]_{\e ^2}
u(m+ \beta _{qn}^{''}) \nonumber
\end{eqnarray}
where 
\begin{eqnarray}
&&\al _j:= -\e_{jj} + \sum _{k=1}^{j-1}
(-\e _{jk} +\e _{j+1,k} -\e _{kj} +\e _{k,j-1}),
\q (2 \leq j \leq n-1), \nonumber \\
&& \al _n:= -\e_{nn} +2 \sum _{k=1}^{n-1}
(-\e _{nk} +\e _{k,n-1}), \nonumber\\
&&\beta _{qj}:= -\e_{qj} + \sum _{p=0}^{q-1}
(-\e _{jp} +\e _{j+1,p} -\e _{pj} +\e _{p,j-1}),
\q (1  \leq q<j \leq n-1), \nonumber \\
&&\beta _{qn}:= -\e_{qn} +2 \sum _{p=0}^{q-1}
(-\e _{np} +\e _{p,n-1}), \q (1 \leq q <n) \nonumber\\
&&\beta _{qj} ^{'}:= \e_{q,j-1} -\e _{jq} -\e _{qj} + \sum _{p=0}^{q-1}
(-\e _{jp} +\e _{j+1,p} -\e _{pj} +\e _{p,j-1}),\nonumber \\
&&\q \qq  (1 \leq q<j \leq n-1), \nonumber \\
&&\beta _{qn} ^{'}:= \e_{q,n-1} -\e _{nq} -\e _{qn} +2 \sum _{p=0}^{q-1}
(-\e _{np} +\e _{p,n-1}), \q (1 \leq q<n),\nonumber \\
&&\beta _{qn} ^{''}:=2 \e_{q,n-1} -2\e _{nq} -\e _{qn} +2 \sum _{p=0}^{q-1}
(-\e _{np} +\e _{p,n-1}), \q (1 \leq q<n). \nonumber
\end{eqnarray}
\end{lem}
{\sl Proof.}
We prove that
for any $m=(m_{ij}) \in M,\, i,j \in I, d \in \z$, 
\begin{equation}
x_{ij}.u(m)=u(m-\e _{ij}), \q \{z_{ij}\}_{\e^d} u(m)=
[d^{-1}(m_{ij}+b_{ij})]_{\e ^d}u(m). 
\label{Lem2.2A} 
\end{equation}
By (\ref{Prop2.1X1}),(\ref{Prop2.1X2}),(\ref{Prop2.1XZ3}), 
we have
\begin{eqnarray}
x_{ij}.u(m)&=&\e^{a_{ij}}X_{ij}(u_{m_{11}}^{(11)} \otimes
\cdots \otimes u_{m_{ij}}^{(ij)} \otimes \cdots \otimes
u_{m_{nn}}^{(nn)})
=u_{m_{11}}^{(11)} \otimes
\cdots \otimes X u_{m_{ij}}^{(ij)} \otimes \cdots \otimes
u_{m_{nn}}^{(nn)} \nonumber \\
&=&u_{m_{11}}^{(11)} \otimes
\cdots \otimes u_{m_{ij}-1}^{(ij)} \otimes \cdots \otimes
u_{m_{nn}}^{(nn)}=u(m-\e_{ij}).
\nonumber
\end{eqnarray} 
Similarly, by (\ref{Prop2.1Z1}),(\ref{Prop2.1Z2}),(\ref{Prop2.1XZ3}), 
\begin{eqnarray}
z_{ij}.u(m)&=&\e^{b_{ij}}z_{ij}.(u_{m_{11}}^{(11)} \otimes
\cdots \otimes u_{m_{ij}}^{(ij)} \otimes \cdots \otimes
u_{m_{nn}}^{(nn)})
\nonumber \\
&=&\e^{b_{ij}}u_{m_{11}}^{(11)} \otimes
\cdots \otimes Zu_{m_{ij}}^{(ij)} \otimes \cdots \otimes
u_{m_{nn}}^{(nn)}
\nonumber \\
&=&\e^{b_{ij}}u_{m_{11}}^{(11)} \otimes
\cdots \otimes \e ^{m_{ij}}u_{m_{ij}}^{(ij)} \otimes \cdots \otimes
u_{m_{nn}}^{(nn)} =\e ^{b_{ij}+m_{ij}}u(m).
\nonumber 
\end{eqnarray}
Therefore , 
\begin{eqnarray}
\{z_{ij}\}_{\e ^d}u(m)&=& \displaystyle \frac{z_{ij}-z_{ij}^{-1}}
{\e ^d -\e ^{-d}}u(m)= 
 \displaystyle \frac{1}{\e ^d -\e ^{-d}}(z_{ij}u(m)-z_{ij}^{-1}u(m)) 
\nonumber \\
&=&\displaystyle \frac{\e ^{b_{ij}+m_{ij}}- \e ^{-b_{ij}-m_{ij}}}
{\e ^d -\e ^{-d}}u(m)=[d^{-1}(m_{ij}+b_{ij})]_{\e ^d} u(m). 
\nonumber
\end{eqnarray} 
We calculate the actions of $e_i$ on 
$u(m)$ by using (\ref{Lem2.2A}) 
and Theorem \ref{Theorem 2.1A}.\\
Case 1) $j=1$:We have 
\begin{eqnarray}
e_1 u&=&\sum _{m \in M} c_m F_{11} u(m)
=\sum _{m \in M} c_m \{z_{11}^{-1}\}x_{11} u(m) \nonumber\\
&=&\sum _{m \in M} c_m \{z_{11}^{-1}\} u(m-\e_{11})
=\sum _{m \in M} c_m [-(m_{11}-1+b_{11})] u(m-\e_{11}) \no\\
&=&\sum _{m \in M} c_m [-(m_{11}-1+1)] u(m-\e_{11}) 
=\sum _{m \in M} c_m [-m_{11}] u(m-\e_{11}). \no
\end{eqnarray}
Case 2) $2 \leq j \leq n-1$: We have 
\begin{eqnarray}
F_{jj}u(m)&=&\{z_{jj}^{-1}\}x_{jj}u(m)=
[-(m_{jj}-1+b_{ij})]u(m-\e _{jj})\no \\
&=&[-(m_{jj}-1+1)]u(m-\e _{jj})=[-m_{jj}]u(m-\e _{jj}). \no 
\end{eqnarray}
Therefore  \\
\begin{eqnarray}
&&(\prod _{k=1}^{j-1}D_{kj})F_{jj}u(m)=[-m_{jj}]
(\prod _{k=1}^{j-1}x_{k,j-1}^{-1}x_{kj}x_{j+1,k}^{-1}x_{jk}
)u(m-\e _{jj}) \no\\
&=&[-m_{jj}]u(m-\e _{jj}+\sum _{k=1}^{j-1}
(-\e _{jk}+\e _{j+1,k}-\e _{kj}+\e _{k,j-1}))
=[-m_{jj}]u(m+\al _j). \no
\end{eqnarray}
\q On the other hand, for any $q (1 \leq q \leq j-1)$, we have
\begin{eqnarray}
C_{qj}u(m)
&=&\{z_{qj}^{-1}z_{q,j-1}\}x_{qj}u(m)+\{z_{jq}^{-1}z_{j+1,q}\}
x_{qj}x_{jq}x_{q,j-1}^{-1}u(m)\no \\
&=&\{z_{qj}^{-1}z_{q,j-1}\}u(m-\e _{qj})+\{z_{jq}^{-1}z_{j+1,q}\}
u(m+\e _{q,j-1}-\e _{jq}-\e _{qj}) \no\\
&=&[-(m_{qj}-1+b_{qj})+(m_{q,j-1}+b_{q,j-1})]u(m-\e _{qj})\no \\
&&+[-(m_{jq}-1+b_{jq})+(m_{j+1,q}+b_{j+1,q})]
u(m+\e _{q,j-1}-\e _{jq}-\e _{qj}) \no\\
&=&[m_{q,j-1}-m_{qj}]u(m-\e _{qj}) +[m_{j+1,q}-m_{jq}]
u(m+\e _{q,j-1}-\e _{jq}-\e _{qj}), \no
\end{eqnarray}
where 
the last equality is due to $b_{jq}=b_{qj}+1=b_{q,j-1}+2$. 
Thus, we obtain
\begin{eqnarray}
&&(\prod _{k=1}^{q-1}D_{pj})C_{qj}u(m)
=[m_{q,j-1}-m_{qj}](\prod_{p=0}^{q-1}x_{p,j-1}^{-1}x_{pj}
x_{j+1,p}^{-1}x_{jp})u(m-\e _{qj}) \no\\
&& +[m_{j+1,q}-m_{jq}]
(\prod_{p=0}^{q-1}x_{p,j-1}^{-1}x_{pj}x_{j+1,p}^{-1}x_{jp})
u(m+\e _{q,j-1}-\e_{jq}-\e_{qj})\no \\
&=&[m_{q,j-1}-m_{qj}]u(m-\e _{qj}+
\sum _{p=0}^{q-1}(-\e _{jp}+\e _{j+1,p}-\e _{pj}+\e _{p,j-1})) \no\\
&& +[m_{q,j-1}-m_{qj}]u(m-\e _{qj}-\e_{jq}+\e _{q,j-1}+
\sum _{p=0}^{q-1}(-\e _{jp}+\e _{j+1,p}-\e _{pj}+\e _{p,j-1})) \no\\
&=&[m_{q,j-1}-m_{qj}]u(m+\beta _{qj})+
[m_{j+1,q}-m_{jq}]u(m+\beta _{qj}^{'}). \no
\end{eqnarray}
\q Therefore, 
\begin{eqnarray}
e_j . u&=&
\sum _{m \in M} c_m(\prod _{k=1}^{j-1}D_{kj})F_{jj}u(m) 
+\sum _{m \in M} \sum _{q=1}^{j-1}
 c_m(\prod _{p=0}^{q-1}D_{pj})C_{qj}u(m) \no\\
&=&\sum _{m \in M} c_m[-m_{jj}]u(m+\al _j) 
+\sum _{m \in M} \sum _{q=1}^{j-1}
 c_m[m_{q,j-1}-m_{qj}]u(m+\beta _{qj}) \no\\ 
&& +\sum _{m \in M} \sum _{q=1}^{j-1}
 c_m[m_{j+1,q}-m_{jq}]u(m+\beta _{qj}^{'}). \no
\end{eqnarray}
The case $j=n$ is shown by the similar way to Case 2.\qed 

Next, we prove the existence and 
uniqueness of the primitive vector. 

\begin{pro}
\label{Proposition 2.2A}
For $\g=\texttt{sp} (2n,\bbC) \q (n \geq 2),$ and 
$\l = (\l _1 , \cdots , \l_n) \in \bbC ^n$, 
set $a^{(0)}=(\e ^{a_{ij}}),b^{(0)}=(\e ^{b_{ij}})$ as 
Lemma \ref{Lemma 2.2A}, and let $(\Phi_{\l,a^{(0)},b^{(0)}}, V)$ 
be the representation as in Theorem \ref{Theorem 2.1A}.
A vector $u\in V$ 
satisfies that $e _i u=0$ for any $i \in I$
if and only if $u \in \bbC u(0)$,
where $u(0)=u(0, \cdots , 0) \in V$. 
(We call these vectors ``primitive vectors''). 
\end{pro}
{\sl Proof}
``If part'' is obvious by Lemma \ref{Lemma 2.2A}. 
So we prove ``only if part''.
First, we define $\{r_i\}_{1\leq i\leq n^2}(=I\times I)$ 
inductively as follows:
$r_1:=(1,1)$ and if $r_s=(i,j)$, then 
\begin{eqnarray}
&& r_{s+1}:=(i,j+1) \q (1 \leq j <n, \, 1 \leq i \leq j+1), \nonumber\\
&& r_{s+1}:=(n,i) \q (1 \leq i < j=n), \nonumber \\
&& r_{s+1}:=(i-1,j) \q (1 \leq j < n, \, j+1<i \leq n).
\label{Prop2.2A} 
\end{eqnarray}
And we define $M_s:=\{m \in M|m_{r_1}=m_{r_2} \cdots =m_{r_s}=0\}$
$(1 \leq s \leq n^2)$. 
So we have
\[ \{(0)\}=M_{n^2} \subset M_{n^2-1} \subset \cdots \subset 
M_1 \subset M.\]
 Now, assume that 
$e_iu=0$ for any $i \in I$ and 
set $u=\sum _{m \in M}c_m u(m) \in V \q (c_m \in \bbC)$.
We shall prove that 
$u=\sum _{m \in M_s}c_m u(m)$ for any $1 \leq s \leq n^2$
 by induction on $s$. 
Indeed, if we can prove this, then we have 
$u=\sum _{m \in M_{n^2}}c_mu(m)=c_0u(0) \in \bbC u(0)$. \\
\q Since $e_1u=0$, by Lemma \ref{Lemma 2.2A},  
we have $0=\sum _{m \in M}c_m[-m_{11}]u(m-\e _{11})$.
Since the vectors $\{u(m-\e _{11})|m \in M\}$ are linearly independent, 
$c_m [-m_{11}]=0$ for any $m \in M$. 
Therefore if $0\neq m_{11} (=m_{r_1})$, then $c_m =0$.
Hence $u=\sum _{m \in M_1} c_m u(m)$. 
Now we assume that
 $u=\sum _{m \in M_s} c_m u(m)$ for  $1 \leq s <n^2$, 
and $r_s =(i,j)$. \\
\q Case 1) $1 \leq i \leq j \leq n-2$: \\
In this case $r_{s+1} =(i,j+1)$. 
Let $m \in M_s$ then 
\[m_{qj}=m_{q,j+1}=m_{j+2,q}=m_{j+1,q}=0
\q  (1 \leq q \leq i-1).  \]
Since $e_ju=0$, by Lemma \ref{Lemma 2.2A}, we have
\begin{eqnarray}
&&0=e_{j+1}u= \nonumber \\
&&\sum _{m \in M} c_m [-m_{j+1,j+1}]u(m+\al _{j+1})
+\sum _{m \in M}\sum _{q=i}^j c_m [m_{qj}-m_{q,j+1}]
u(m+\beta _{q,j+1}) \nonumber \\
&&\qq +\sum _{m \in M}\sum _{q=i}^j c_m [m_{j+2,q}-m_{j+1,q}]
u(m+\beta _{q,j+1}^{'}). \nonumber
\end{eqnarray}
\q On the other hand, by Lemma \ref{Lemma 2.2A}, for any  $m \in M_s$, 
we have 
\begin{eqnarray*}
&&(m+\al _{j+1})_{ij}=
(m-\e _{j+1,j+1}+\sum _{k=1}^j
(-\e _{j+1,k}-\e_{k,j+1}+\e _{kj}))_{ij} 
=(m+\sum _{k=1}^j \e _{kj})_{ij}\\
&&=(m+\e_{ij})_{ij}
=m_{ij}+1=1,
\end{eqnarray*} 
since $m \in M_s, 0=m_{r_s}=m_{ij}$.
Similarly , for any $q(i \leq q \leq j)$, we have 
\begin{eqnarray*}
(m+\beta _{q,j+1})_{ij}&=&(m-\e _{q,j+1}+\sum _{p=0}^{q-1}
(-\e _{j+1,p}+\e _{j+2,p}-\e _{p,j+1}+\e _{pj}))_{ij}
=(\sum _{p=0}^{q-1}\e _{pj})_{ij}\\
(m+\beta _{q,j+1}^{'})_{ij} 
&=&(m-\e _{q,j+1}-\e _{j+1,q}+\e_{qj}
+\sum _{p=0}^{q-1}
(-\e _{j+1,p}+\e _{j+2,p}-\e _{p,j+1}+\e _{pj}))_{ij}\\
&=&(\sum _{p=0}^{q}\e _{pj})_{ij}.
\end{eqnarray*}
Therefore 
\begin{eqnarray}
&&(m+\al _{j+1})_{ij}=1, 
\qq (m+\beta _{q,j+1})_{ij}=1 \q (i+1 \leq q \leq j), \nonumber \\
&& (m+\beta _{q,j+1})_{ij}=0 \q (q=i),
\qq  (m+\beta _{q,j+1}^{'})_{ij}=1 \q (i \leq q \leq j). \nonumber
\end{eqnarray}
Thus 
\begin{eqnarray*}
&&\sum _{m \in M_s} \bbC u(m+\al _{j+1})
+\sum _{m \in M_s} \sum _{q=i}^{j} \bbC u(m+\beta _{q,j+1}) 
 +\sum _{m \in M_s} \sum _{q=i}^{j} \bbC u(m+\beta _{q,j+1}^{'})\\
&&=\{\bigoplus _{m \in M_s}  \bbC u(m+\beta _{i,j+1}) \}
\bigoplus \{\sum _{m \in M_s} \bbC u(m+\al _{j+1}) \\ 
&&\q+\sum _{m \in M_s} \sum _{q=i+1}^{j} \bbC u(m+\beta _{q,j+1})+
\sum _{m \in M_s} \sum _{q=i}^{j} 
\bbC u(m+\beta _{q,j+1}^{'}) \}.
\end{eqnarray*}
Then since the vectors 
$\{u(m+\beta_{i,j+1}) \, |\, m \in M_s$ \} are 
linearly independent, we have that 
$0=c_m [m_{ij}-m_{i,j+1}]=c_m [-m_{i,j+1}]$ \q 
for any $m \in M_s$, which implies 
 if $0 \neq m_{i,j+1}=m_{r_{s+1}}$ then $c_m=0$ .
Therefore $u=\sum _{m \in M_{s+1}}c_m u(m)$. 

Case 2) \q $1 \leq i \leq j = n-1$: \\
In this case , $r_s=(i,n)$.
Let $m \in M$ then 
\[m_{qn}=m_{nq}=m_{q,n-1}=0 \q (1 \leq q \leq i-1). \]
Thus, since $e_n u=0$, by Lemma \ref{Lemma 2.2A} 
\begin{eqnarray}
&& 0=e_n u \nonumber \\
&&=\sum _{m \in M_s} c_m [-2m_{nn}]_{\e ^2}u(m+\al _n) 
 +\sum _{m \in M_s} \sum _{q=i}^{n-1} 
 c_m [2(m_{qn}-m_{nq})]_{\e ^2}u(m+\beta _{qn}) \nonumber \\
&&  +\sum _{m \in M_s} \sum _{q=i}^{n-1}  
 c_m [m_{q,n-1}-m_{nq}]u(m+\beta _{qn}^{'}) 
  +\sum _{m \in M_s} \sum _{q=i}^{n-1} 
 c_m [2(m_{qn}-m_{nq})]_{\e ^2}u(m+\beta _{qn}^{''}). \nonumber 
\end{eqnarray}
On the other hand, by the similar way to the Case 1, we have
\begin{eqnarray}
&&(m+\al _n)_{i,n-1}=2, 
\qq (m+\beta _{qn})_{i,n-1}=2 \q (i+1 \leq q \leq n-1), \nonumber \\
&& (m+\beta _{qn})_{i,n-1}=0 \q (q=i),
\qq (m+\beta _{qn}^{'})_{i,n-1}=2 \q (i+1 \leq q \leq n-1), \nonumber \\ 
&&(m+\beta _{qn}^{'})_{i,n-1}=1 \q (q=i),
\qq (m+\beta _{qn}^{''})_{i,n-1}=2 \q (i \leq q \leq n-1).\nonumber
\end{eqnarray}
Then, we have 
\begin{eqnarray*}
&&\sum _{m \in M_s} \bbC u(m+\al _{n})
+\sum _{m \in M_s} \sum _{q=i}^{n-1} \bbC u(m+\beta _{qn})
+\sum _{m \in M_s} \sum _{q=i}^{n-1} \bbC u(m+\beta _{qn}^{'})
+\sum _{m \in M_s} \sum _{q=i}^{n-1} \bbC u(m+\beta _{qn}^{''}) \\
&&=\{\bigoplus _{m \in M_s}  \bbC u(m+\beta _{in}) \}
\bigoplus \{\bigoplus _{m \in M_s}  \bbC u(m+\beta _{in}^{'}) \} 
\bigoplus \{\sum _{m \in M_s} \bbC u(m+\al_{n})
+\sum _{m \in M_s} \sum _{q=i+1}^{n-1} \bbC u(m+\beta _{qn}) \\
&&+ \sum _{m \in M_s} \sum _{q=i+1}^{n-1} \bbC u(m+\beta _{qn}^{'})
+\sum _{m \in M_s} \sum _{q=i}^{n-1} \bbC u(m+\beta _{qn}^{''}) \}. 
\end{eqnarray*}
Hence by the linearly independence of the vectors 
$\{u(m+\beta_{i,n}) \, | \, m \in M_s\}$,
$\{u(m+\beta_{i,n}^{'}) \, | \, m \in M_s\}$, we have that 
$0=c_m[2(m_{in}-m_{ni})]=c_m[m_{i,n-1}-m_{ni}]=c_m[-m_{ni}]$ 
for any $m \in M_s$. 
Thus if $m_{in} \neq 0$ or $m_{ni} \neq 0 $,
then $c_m =0$ \q  for any $m \in M_s$. 
Therefore $u= \sum _{m \in M_{s+2}}c_mu(m)
= \sum _{m \in M_{s+1}}c_mu(m)$. 

Case 3) $1 \leq i \leq j=n$.
This case is shown by the similar way to Case2.

Case 4) $1 \leq j < i$ and $j<i-1$: \\
In this case $r_{s+1} =(i-1,j)$. 
For $m \in M_s$, we have \\
\[m_{q,i-2}=m_{q,i-2}=0 \q  (1 \leq q \leq j),
\qq m_{iq}=m_{i-1,q}=0 \q  (1 \leq q \leq j-1). \]
Since $e_{i-1}u=0$, by Lemma \ref{Lemma 2.2A}, 
\begin{eqnarray}
&&0=e_{i-1}u= \nonumber \\
&&\sum _{m \in M_s}c_m [-m_{i-1,i-1}]u(m+\al _{i-1}) 
 + \sum _{m \in M_s} \sum _{q=j+1}^{i-2}
c_m [m_{q,i-2}-m_{q,i-1}]u(m+\beta _{q,i-1})\nonumber \\
&& +  \sum _{m \in M_s} \sum _{q=j}^{i-2}
c_m [m_{iq}-m_{i-1,q}]u(m+\beta _{q,i-1}^{'}) \nonumber
\end{eqnarray}
On the other hand,
\begin{eqnarray} 
&& (m+\al _{i-1})_{ij}=1, 
\qq (m+\beta _{q,i-1})_{ij}=1 \q  (j+1 \leq q \leq i-2), \nonumber \\
&&(m+\beta _{q,i-1}^{'})_{ij}=1 \q  (j+1 \leq q \leq i-2),
\qq (m+\beta _{q,i-1}^{'})_{ij}=0 \q  (q=j).\nonumber
\end{eqnarray}
\q Thus, $0=c_m [m_{ij}-m_{i-1,j}]=c_m [-m_{i-1,j}]=0$  
\q for any $m \in M_s$. 
Therefore if $0 \neq m_{i-1,j}=m_{r_{s+1}}$ then $c_m=0$
\q for any $m \in M_s$. 
So $u= \sum _{m \in M_s}c_m u(m)$. \\
\q Case 5) $1 \leq j<i \leq n$ and $j=i-1$: \\
In this case $r_{s+1} =(i,i)$. For $m \in M_s$, we have
\[m_{q,i-1}=m_{qi}=m_{iq}=m_{i+1,q}=0 \q (1 \leq q \leq i-1)\]
Since $e_i u=0$, by Lemma \ref{Lemma 2.2A}, we get
\begin{eqnarray}
&&0=e_i u=\sum _{m \in M_s}c_m [-m_{ii}]
u(m+\al _i) \q (i \neq n) \nonumber \\
&&0=e_n u=\sum _{m \in M_s}c_m [-2m_{nn}]_{\e ^2}
u(m+\al _n) \q (i = n) \nonumber
\end{eqnarray}
\q Hence, if $m_{ii} \neq 0$ then $c_m=0$ \q for any $m \in M_s$. 
So $u=\sum  _{m \in M_{s+1}}c_m u(m)$. \qed 

\subsubsection{$\texttt{so}(2n+1,\bbC)$-case}
By the similar way to the proof of 
Lemma \ref{Lemma 2.2A} and 
Proposition \ref{Proposition 2.2A}, 
we can prove the following lemma and proposition. 
\begin{lem}
\label{Lemma 2.2B}
For $\g :=\texttt{so} (2n+1,\bbC)  \q (n \geq 3),
\l :=(\l _1 , \cdots , \l _n) \in \bbC ^n$,
and any $i,j \in I$, set 
\begin{eqnarray}
&& a_{ij}:=0,
\qq  b_{ij}:=1-i+j \q (i \leq j \leq n-1 ), \nonumber \\
&& b_{in}:=2n+1-2i, 
\qq  b_{ij}:=2n+1-i-j \q (i>j), \nonumber
\end{eqnarray}
and $a^{(0)}:=(\e ^{a_{ij}})_{i,j=1}^n,
b^{(0)}:=(\e ^{b_{ij}})_{i,j=1}^n \in (\bbC ^{\times})^{n^2}$.
Let $(\Phi _{\l , a^{(0)} , b^{(0)}}, V)$ 
be the representation of Theorem \ref{Theorem 2.1B}. \\
For $u= \sum _{m \in M} c_m u(m)\in V \q (c_m \in \bbC)$, we have 
\begin{eqnarray}
&&e_1 . u= \sum _{m \in M}c_m [-2m_{11}]_{\e^2}u(m-\e _{11}),\nonumber \\
&&e_j . u= \sum _{m \in M}c_m [-2m_{jj}]_{\e^2}u(m+ \al _j)  
+ \sum _{m \in M} \sum_{q=1}^{j-1}c_m [2(m_{q,j-1}-m_{qj})]_{\e^2}
u(m+ \beta _{qj}) \nonumber \\
&&\qq \q + \sum _{m \in M} \sum_{q=1}^{j-1}c_m [2(m_{j+1,q}-m_{jq})]_{\e^2}
u(m+ \beta _{qj}^{'}), 
\q  (2 \leq j \leq n-1), \nonumber \\
&&e_n . u= \sum _{m \in M}c_m [-m_{nn}]u(m+\al_n) 
 + \sum _{m \in M} \sum_{q=1}^{n-1}c_m [2m_{q,n-1}-m_{qn}]
u(m+ \beta _{qn}) \nonumber \\
&&\qq \q  + \sum _{m \in M} \sum_{q=1}^{n-1}c_m [m_{qn}-2m_{nq}]
u(m+ \beta _{qn}^{'}), \nonumber  
\end{eqnarray}
where 
\begin{eqnarray}
&&\al _j:= -\e_{jj} + \sum _{k=1}^{j-1}
(-\e _{jk} +\e _{j+1,k} -\e _{kj} +\e _{k,j-1}),
\q (2 \leq j \leq n-1),  \nonumber \\
&&\al _n:= -\e_{nn} + \sum _{k=1}^{n-1}
(-\e _{nk} +\e _{k,n-1}),\no
\end{eqnarray}
\begin{eqnarray}
&&\beta _{qj}:= -\e_{qj} + \sum _{p=0}^{q-1}
(-\e _{jp} +\e _{j+1,p} -\e _{pj} +\e _{p,j-1}),
\q (1 \leq q<j \leq n-1), \nonumber  \\
&&\beta _{qn}:= -\e_{qn} + \sum _{p=0}^{q-1}
(-\e _{np} +\e _{p,n-1}), \q (1 \leq q<n), \nonumber \\
&&\beta _{qj} ^{'}:= \e_{q,j-1} -\e _{jq} -\e _{qj} + \sum _{p=0}^{q-1}
(-\e _{jp} +\e _{j+1,p} -\e _{pj} +\e _{p,j-1}), 
\q (1 \leq q<j \leq n-1),  \nonumber \\
&&\beta _{qn} ^{'}:= \e_{q,n-1} -\e _{nq} -\e _{qn} + \sum _{p=0}^{q-1}
(-\e _{np} +\e _{p,n-1}), \q (1 \leq q<n).  \nonumber 
\end{eqnarray}
\end{lem}
\begin{pro}
\label{Proposition 2.2B}
For $\g=\texttt{so}(2n+1,\bbC) \, (n \geq 3)$ and 
$\l =(\l_1 , \cdots , \l_n) \in \bbC ^n$,
set $a^{(0)}=(\e^{a_{ij}}), 
b^{(0)}=(\e^{b_{ij}})$ as Lemma \ref{Lemma 2.2B}, 
and let $(\Phi _{\l , a^{(0)} , b^{(0)}}, V)$ 
be the representation as in Theorem \ref{Theorem 2.1B}. 
A vector $u=\sum _{m \in M} c_m u(m) \in V \, (c_m \in \bbC)$
satisfies the condition 
$e_i u =0$ for any $i \in I$ if and only if 
$u \in \bbC u(0)$. 
\end{pro}
{\sl Sketch of the proof.}
``If part'' is obvious from Lemma \ref{Lemma 2.2B}. 
So we prove ``only part''.
We define $\{r_s\}_{s=1}^{n^2}$ by the same way 
as (\ref{Prop2.2A}) 
and $M_s=\{m \in M | m_{r_1} = \cdots =m_{r_s}=0 \}$ 
for any $s(1 \leq s \leq n^2)$. 
Now we prove that $u=\sum_{m \in M_{s}} c_m u(m)$ for any 
$s (1 \leq s \leq n^2)$ by induction on $s$. 
By using $e_1 u=0$, we can prove $u=\sum_{m \in M_{1}} c_m u(m)$. 
Next, we assume that 
$u=\sum_{m \in M_{s}}c_m u(m)$ for a $s (1 \leq s<n^2)$ and $r_s=(i,j)$. 
Then we prove by the similar way to the proof of 
Proposition \ref{Proposition 2.2A}. \\
\q Case 1) $1 \leq i \leq j \leq n-1$ :  \\
We use $e_{j+1}u=0$ and for any $m \in M_s$, 
\begin{eqnarray}
&& (m+\al _{j+1})_{ij}=1, 
\qq (m+\beta _{q,j+1})_{ij}=1 \q  (q>i), \nonumber \\
&&(m+\beta _{q,j+1})_{ij}=0 \q  (q=i),
\qq(m+\beta _{q,j+1}^{'})_{ij}=1 \q  (q \geq i). \nonumber
\end{eqnarray}
\q Case 2) $1 \leq i \leq j=n$ :  \\
We use $e_{n}u=0$ and for any $m \in M_s$, 
\begin{eqnarray}
&& (m+\al _{n})_{in}=0, 
\qq(m+\beta _{qn})_{in}=0 \q  (q>i),\nonumber \\
&&(m+\beta _{qn}^{'})_{in}=0 \q  (q>i),
\qq (m+\beta _{qn}^{'})_{in}=l-1 \q  (q=  i). \nonumber
\end{eqnarray}
\q Case 3) $1 \leq j<i-1 < n$ :  \\
We use $e_{j+1}u=0$ and for any $m \in M_s$, 
\begin{eqnarray}
&& (m+\al _{i-1})_{ij}=1, 
\qq (m+\beta _{q,i-1})_{ij}=1 \q  (q \geq j+1), \nonumber \\
&& (m+\beta _{q,i-1}^{'})_{ij}=1 \q  (q \geq j+1),
\qq (m+\beta _{q,i-1}^{'})_{ij}=0 \q  (q =j). \nonumber
\end{eqnarray}
\q Case 4) $1 \leq j=i-1 < n$ : We use $e_{j+1}u=0$.  \qed

\subsubsection{$\texttt{so}(2n, \bbC)$-case}
By the similar way to the proof of 
Lemma \ref{Lemma 2.2A} and 
Proposition \ref{Proposition 2.2A}, 
we can prove the following lemma and proposition. 
\begin{lem}
\label{Lemma 2.2C}
For $\g :=\texttt{so} (2n,\bbC)  \, (n \geq 4),
\l :=(\l _1 , \cdots , \l _n) \in \bbC ^n$,
and  any $i,j \in I$, set 
\begin{eqnarray}
 &&a_{ij}:=0,
\qq  b_{ij}:=1-i+j \q (i \leq j \leq n-1 ), \nonumber \\
&& b_{in}:=n-i, 
\qq  b_{ij}:=2n-i-j \q (j<i \leq n-1),\nonumber
\end{eqnarray}
and $a^{(0)}:=(\e ^{a_{ij}})_{1 \leq i \leq n-1, 1 \leq j \leq n},
b^{(0)}:=(\e ^{b_{ij}})_{1 \leq i \leq n-1, 1 \leq j \leq n}
 \in (\bbC ^{\times})^{n(n-1)}$.
Let $(\Phi _{\l , a^{(0)} , b^{(0)}}, V)$ 
be the representation as in Theorem \ref{Theorem 2.1C}. 
For  $u= \sum _{m \in M} c_m u(m)\in V \, (c_m \in \bbC)$, we have 
\begin{eqnarray}
&&e_1 . u= \sum _{m \in M}c_m [-m_{11}]_{\e^2}u(m-\e _{11}), 
\nonumber\\
&&e_j . u= \sum _{m \in M}c_m [-m_{jj}]_{\e^2}u(m+ \al _j) 
+ \sum _{m \in M} \sum_{q=1}^{j-1}c_m [m_{q,j-1}-m_{qj}]
u(m+ \beta _{qj}) \nonumber\\
&&\qq \q + \sum _{m \in M} \sum_{q=1}^{j-1}c_m [m_{j+1,q}-m_{jq}]
u(m+ \beta _{qj}^{'}),  
\q (2 \leq j \leq n-1), \nonumber
\end{eqnarray}
\begin{eqnarray}
&&e_{n-1} . u= \sum _{m \in M}c_m [-m_{n-1,n-1}]u(m+ \al _{n-1})  \nonumber\\
&&+ \sum _{m \in M} \sum_{q=1}^{n/2-1}
\{c_m [m_{2q-1,n-2}-m_{2q-1,n-1}]u(m+ \beta _{q,n-1})
+c_m [m_{2q-1,n}-m_{n-1,2q-1}]u(m+ \beta _{q,n-1}^{'}) \nonumber \\
&&+c_m [m_{2q,n-2}-m_{2q,n}]u(m+ \gamma _{q,n-1}) 
+c_m [m_{2q,n-1}-m_{n-1,2q}]u(m+ \gamma _{q,n-1}^{'})\}, 
\q (n; even)\no
\end{eqnarray}
\begin{eqnarray}
&&e_{n} . u= \sum _{m \in M}c_m [-m_{n-1,n}]u(m+ \al _{n}) \nonumber\\ 
&&+ \sum _{m \in M} \sum_{q=1}^{n/2-1}
\{c_m [m_{2q,n-2}-m_{2q,n-1}]u(m+ \beta _{q,n})
+c_m [m_{2q,n}-m_{n-1,2q}]u(m+ \beta _{q,n}^{'}) \nonumber \\
&&+ c_m [m_{2q-1,n-2}-m_{2q-1,n}]u(m+ \gamma _{q,n}) 
 + c_m [m_{2q-1,n-1}-m_{n-1,2q-1}]u(m+ \gamma _{q,n}^{'})\}, 
\q (n; even)\nonumber 
\end{eqnarray}
\begin{eqnarray}
&&e_{n-1} . u= \sum _{m \in M}c_m [-m_{n-1,n}]u(m+ \al _{n-1})\nonumber \\ 
&&+ \sum _{m \in M} \sum_{q=1}^{(n-1)/2-1}
c_m [m_{2q-1,n-2}-m_{2q-1,n-1}]
u(m+ \beta _{q,n-1}) \nonumber\\
&& + \sum _{m \in M} \sum_{q=1}^{(n-1)/2}c_m [m_{2q-1,n}-m_{n-1,2q-1}]
u(m+ \beta _{q,n-1}^{'})  \nonumber\\
&&+ \sum _{m \in M} \sum_{q=1}^{(n-3)/2}\{c_m [m_{2q,n-2}-m_{2q,n}]
u(m+ \gamma _{q,n-1}) 
+c_m [m_{2q,n-1}-m_{n-1,2q}]u(m+ \gamma _{q,n-1}^{'})\}, 
\q (n; odd)\nonumber 
\end{eqnarray}
\begin{eqnarray}
&&e_{n} . u= \sum _{m \in M}c_m [-m_{n-1,n-1}]u(m+ \al _{n}) \nonumber\\ 
&&+ \sum _{m \in M} \sum_{q=1}^{(n-3)/2}
\{c_m [m_{2q,n-2}-m_{2q,n-1}]u(m+ \beta _{q,n}) 
+ c_m [m_{2q,n}-m_{n-1,2q}]u(m+ \beta _{q,n}^{'})  \nonumber\\
&&+ \sum _{m \in M} \sum_{q=1}^{(n-1)/2}
\{c_m [m_{2q-1,n-2}-m_{2q-1,n}]u(m+ \gamma _{q,n})
c_m [m_{2q-1,n-1}-m_{n-1,2q-1}]u(m+ \gamma _{q,n}^{'})\}, 
\q (n; odd) \nonumber
\end{eqnarray}
where, 
\begin{eqnarray}
&&\al_j=-\e_{jj}+\sum_{k=1}^{j-1}(-\e_{jk}+\e_{j+1,k}-\e_{kj}+\e_{k,j-1}),
\q  (2 \leq j \leq n-2), \nonumber\\
&&\al_{n-1}=-\e_{n-1,n-1}+\sum_{k=1}^{n/2-1}
(-\e_{n-1,2k-1}+\e_{2k-1,n}-\e_{2k-1,n-1}+\e_{2k-1,n-2})\nonumber \\
&&\q \qq \qq +\sum_{k^{'}=1}^{n/2-1}
(-\e_{n-1,2k^{'}}+\e_{2k^{'},n-1}-\e_{2k^{'},n}+\e_{2k^{'},n-2}),
\q (n; even)\nonumber \\
&&\q \al_{n}=-\e_{n-1,n}+\sum_{k=1}^{n/2-1}
(-\e_{n-1,2k-1}+\e_{2k-1,n-1}-\e_{2k-1,n}+\e_{2k-1,n-2})\nonumber \\
&&\q \qq \qq +\sum_{k^{'}=1}^{n/2-1}
(-\e_{n-1,2k^{'}}+\e_{2k^{'},n-1}-\e_{2k^{'},n}+\e_{2k^{'},n-2}),
\q (n; even)\nonumber 
\end{eqnarray}
\begin{eqnarray}
&&\q \al_{n-1}=-\e_{n-1,n}+\sum_{k=1}^{(n-3)/2}
(-\e_{n-1,2k}+\e_{2k,n-1}-\e_{2k,n}+\e_{2k,n-2}) \nonumber\\
&&\q \qq \qq +\sum_{k^{'}=1}^{(n-1)/2}
(-\e_{n-1,2k^{'}-1}+\e_{2k^{'}-1,n}-\e_{2k^{'}-1,n-1}+\e_{2k^{'}-1,n-2}),
\q (n; odd)\nonumber \\
&&\al_{n}=-\e_{n-1,n-1}+\sum_{k=1}^{(n-1)/2}
(-\e_{n-1,2k-1}+\e_{2k-1,n-1}-\e_{2k-1,n}+\e_{2k-1,n-2}) \nonumber\\
&&\q \qq \qq +\sum_{k^{'}=1}^{(n-3)/2}
(-\e_{n-1,2k^{'}}+\e_{2k^{'},n}-\e_{2k^{'},n-1}+\e_{2k^{'},n-2}),
\q (n; odd)\nonumber 
\end{eqnarray}
\begin{eqnarray}
&&\beta_{qj}=-\e_{qj}+\sum_{p=0}^{q-1}
(-\e_{jp}+\e_{j+1,p}-\e_{pj}+\e_{p,j-1}), 
\q (1 \leq q<j \leq n-2),\nonumber \\
&&\beta_{q,n-1}=-\e_{2q-1,n-1}+\sum_{p^{'}=1}^{q-1}
(-\e_{n-1,2p^{'}}+\e_{2p^{'},n-1}-\e_{2p^{'},n}+\e_{2p^{'},n-2})\nonumber \\ 
&&\q  \qq \qq +\sum_{p=1}^{q-1}
(-\e_{n-1,2p-1}+\e_{2p-1,n}-\e_{2p-1,n-1}+\e_{2p-1,n-2}),\nonumber\\
&&\qq \q (n; even, 1 \leq q < n/2-1), 
 \q (n; odd, 1 \leq q < (n-1)/2-1),\nonumber \\
&&\beta_{qn}=-\e_{2q,n-1}+\sum_{p^{'}=1}^{q}
(-\e_{n-1,2p^{'}-1}+\e_{2p^{'}-1,n-1}
-\e_{2p^{'}-1,n}+\e_{2p^{'}-1,n-2})
\nonumber \\ 
&&\q  \qq \qq +\sum_{p=1}^{q-1}
(-\e_{n-1,2p}+\e_{2p,n}-\e_{2p,n-1}+\e_{2p,n-2}),\nonumber \\
&&\qq \q (n; even, 1 \leq q < n/2-1), 
 \q (n; odd, 1 \leq q < (n-3)/2-1), \nonumber \\
&&\q \beta^{'}_{qj}=-\e_{qj}-\e_{jq}+\e_{q,j-1}+\sum_{p=0}^{q-1}
(-\e_{jp}+\e_{j+1,p}-\e_{pj}+\e_{p,j-1}),  
 \q (1 \leq q<j \leq n-2), \nonumber\\
&&\q \beta_{q,n-1}^{'}=-\e_{2q-1,n-1}-\e_{n-1,2q-1}+\e_{2q-1,n-2}
 +\sum_{p^{'}=0}^{q-1}
(-\e_{n-1,2p^{'}}+\e_{2p^{'},n-1}-\e_{2p^{'},n}+\e_{2p^{'},n-2})\nonumber \\ 
&&\q  \qq \qq +\sum_{p=0}^{q-1}
(-\e_{n-1,2p-1}+\e_{2p-1,n}-\e_{2p-1,n-1}+\e_{2p-1,n-2}),\nonumber\\  
&&\qq \q (n; even, 1 \leq q < n/2-1), 
 \q (n; odd, 1 \leq q < (n-1)/2-1), \nonumber\\ 
&&\q \beta_{qn}^{'}=-\e_{2q,n-1}-\e_{n-1,2q}+\e_{2q,n-2}
+\sum_{p^{'}=0}^{q}(-\e_{n-1,2p^{'}-1}+\e_{2p^{'}-1,n-1}-\e_{2p^{'}-1,n}
+\e_{2p^{'}-1,n-2})\nonumber  \\ 
&&\q  \qq \qq +\sum_{p=0}^{q-1}
(-\e_{n-1,2p}+\e_{2p,n}-\e_{2p,n-1}+\e_{2p,n-2}), \nonumber \\ 
&&\qq \q (n; even, 1 \leq q < n/2-1), 
 \q (n; odd, 1 \leq q < (n-3)/2-1),\nonumber  
\end{eqnarray}
\begin{eqnarray}
&&\q \gamma_{q,n-1}=-\e_{2q,n}+\sum_{p^{'}=1}^{q-1}
(-\e_{n-1,2p^{'}}+\e_{2p^{'},n-1}-\e_{2p^{'},n}+\e_{2p^{'},n-2}) \nonumber \\ 
&&\q  \qq \qq +\sum_{p=1}^{q-1}
(-\e_{n-1,2p-1}+\e_{2p-1,n}
-\e_{2p-1,n-1}+\e_{2p-1,n-2}),\nonumber  \\ 
&&\qq \q (n; even, 1 \leq q < n/2-1), 
 \q (n; odd, 1 \leq q < (n-3)/2-1), \nonumber \\ 
&&\gamma_{qn}=-\e_{2q-1,n-1}+\sum_{p^{'}=1}^{q-1}
(-\e_{n-1,2p^{'}-1}+\e_{2p^{'}-1,n-1}
-\e_{2p^{'}-1,n}+\e_{2p^{'}-1,n-2})
 \nonumber \\ 
&&\q  \qq \qq +\sum_{p=1}^{q-1}
(-\e_{n-1,2p}+\e_{2p,n}-\e_{2p,n-1}+\e_{2p,n-2}), \nonumber \\ 
&&\qq \q (n; even, 1 \leq q < n/2-1), 
 \q (n; odd, 1 \leq q < (n-1)/2-1), \nonumber 
\end{eqnarray}
\begin{eqnarray}
&&\gamma_{q,n-1}^{'}=-\e_{2q,n-1}-\e_{n-1,2q}+\e_{2q,n-2}+\sum_{p^{'}=1}^{q-1}
(-\e_{n-1,2p^{'}}+\e_{2p^{'},n-1}-\e_{2p^{'},n}+\e_{2p^{'},n-2}) \nonumber \\ 
&&\q  \qq \qq +\sum_{p=1}^{q}
(-\e_{n-1,2p-1}+\e_{2p-1,n}-\e_{2p-1,n-1}+\e_{2p-1,n-2}), \nonumber \\ 
&&\qq \q (n; even, 1 \leq q < n/2-1), 
 \q (n; odd, 1 \leq q < (n-3)/2-1), \nonumber \\
&&\gamma_{qn}^{'}=-\e_{2q-1,n}-\e_{n-1,2q-1}+\e_{2q-1,n-2} 
+\sum_{p^{'}=0}^{q-1}
(-\e_{n-1,2p^{'}-1}+\e_{2p^{'}-1,n-1}-\e_{2p^{'}-1,n}+\e_{2p^{'}-1,n-2})
\nonumber  \\ 
&&  \qq \qq +\sum_{p=0}^{q-1}
(-\e_{n-1,2p}+\e_{2p,n}-\e_{2p,n-1}+\e_{2p,n-2}),\nonumber  \\ 
&&\qq \q (n; even, 1 \leq q < n/2-1),  
 \q (n; odd, 1 \leq q < (n-1)/2-1).\nonumber 
\end{eqnarray}
\end{lem}
\begin{pro}
\label{Proposition 2.2C}
For $\g=\texttt{so}(2n,\bbC) \, (n \geq 4)$ and 
$\l =(\l_1 , \cdots , \l_n) \in \bbC ^n$,
set $a^{(0)}=(\e^{a_{ij}}),b^{0}
=(\e^{b_{ij}})$ as Lemma \ref{Lemma 2.2C}, 
and $(\Phi _{\l , a^{(0)} , b^{(0)}}, V)$ 
be the representation as in Theorem \ref{Theorem 2.1C}. 
A vector $u=\sum _{m \in M} c_m u(m) \in V \, (c_m \in \bbC)$
satisfies the condition 
$e_i u =0$ for any $i \in I$ if and only if 
$u \in \bbC u(0)$. 
\end{pro}
{\sl Sketch of the proof.} 
``If part'' is obvious from Lemma \ref{Lemma 2.2C}. 
So we prove ``only part''.
We define $\{r_i\}_{i=1}^{n(n-1)}$
 by the similar way to the previous cases:
$r_1:=(1,1)$ and if $r_s=(i,j)$, then   
\begin{eqnarray}
&&r_{s+1}:=(i,j+1) 
\q (1 \leq i,j \leq n-1, \, 1 \leq i \leq j+1), \nonumber \\
&&r_{s+1}:=(n,i)
\q (1 \leq i \leq n-1, \,  j=n), \nonumber \\
&&r_{s+1}:=(i-1,j) 
\q (1 \leq j <n-1, \,  j+1 \leq i \leq n-1). \nonumber
\end{eqnarray}
We set $M_s=\{m \in M | m_{r_1} = \cdots =m_{r_s}=0 \}$ 
for any $s(1 \leq s \leq n(n-1))$. \\
\q Now we prove that $u=\sum_{m \in M_{s}} c_m u(m)$ for any $1 \leq s
\leq n(n-1)$ by induction on $s$. 
By using $e_1 u=0$, we can prove $u=\sum_{m \in M_{1}} c_m u(m)$. 
Next, we assume that $u=\sum_{m \in M_{s}}c_m u(m)$
for a $s (1 \leq s<n(n-1))$.  
There exist a pair $(i,j)$ such that 
$ (i,j) \neq (n-1,n)$ and $r_s^m=m_{ij}$ for any $m \in M_s$. 
Then we prove by the similar manner to the proof of Proposition 
\ref{Proposition 2.2A}. \\
\q Case 1) $1 \leq i \leq j \leq n-3$:  \\
We use $e_{j+1}u=0$ and for any $m \in M_s$,
\begin{eqnarray} 
&&(m+\al _{j+1})_{ij}=1,
\qq (m+\beta _{q,j+1})_{ij}=1 \q  (q>i), \nonumber \\
&&(m+\beta _{q,j+1})_{ij}=0 \q  (q=i),
\qq(m+\beta _{q,j+1}^{'})_{ij}=1 \q  (q \geq i). \nonumber 
\end{eqnarray}
\q Case 2.1 \, $1 \leq i \leq j \leq n-2$ and $i=2i^{'}+1$:  \\
We use $e_{n-1}u=0$ and for any $m \in M_s$, 
\begin{eqnarray} 
&&(m+\al _{n-1})_{i,n-2}=1, 
\qq (m+\beta _{q,n-1})_{i,n-2}=1 \q  (q>i^{'}),\nonumber \\
&& (m+\beta _{q,n-1})_{i,n-2}=0 \q  (q=i^{'}), 
\qq (m+\beta _{q,n-1}^{'})_{i,n-2}=1 \q  (q \geq i^{'}),\nonumber \\
&&q (m+\gamma _{q,n-1})_{i,n-2}=1 \q  (q \geq i^{'}),
\qq (m+\gamma _{q,n-1}^{'})_{i,n-2}=1 \q  (q \geq i^{'}+1), \nonumber
\end{eqnarray}
and we use $e_{n}u=0$ and for any $m \in M_s$, 
\begin{eqnarray} 
&& (m+\al _{n})_{i,n-2}=1,
\qq (m+\beta _{qn})_{i,n-2}=1 \q  (q \geq i^{'}), \nonumber \\
&&(m+\beta _{qn}^{'})_{i,n-2}=1 \q  (q \geq i^{'}),
\qq(m+\gamma _{qn})_{i,n-2}=1 \q  (q \geq i^{'}+1), \nonumber \\
&&(m+\gamma _{qn})_{i,n-2}=0 \q  (q =i^{'}),
\qq(m+\gamma _{qn}^{'})_{i,n-2}=1 \q  (q \geq i^{'}). \nonumber
\end{eqnarray}
\q Case 2.2) $1 \leq i \leq j \leq n-2$ and $i=2i^{'}$:  \\
We use $e_{n-1}u=0$ and for any $m \in M_s$, 
\begin{eqnarray} 
&& (m+\al _{n-1})_{i,n-2}=1, 
\qq (m+\beta _{q,n-1})_{i,n-2}=1 \q  (q>i^{'}+1), \nonumber \\
&&(m+\beta _{q,n-1}^{'})_{i,n-2}=1 \q  (q \geq i^{'}+1),
\qq (m+\gamma _{q,n-1})_{i,n-2}=1 \q  (q \geq i^{'}+1), \nonumber \\
&&(m+\gamma _{q,n-1})_{i,n-2}=0 \q  (q=  i^{'}),
\qq (m+\gamma _{q,n-1}^{'})_{i,n-2}=1 \q  (q \geq i^{'}), \nonumber
\end{eqnarray}
and we use $e_{n}u=0$ and for any $m \in M_s$,
\begin{eqnarray}  
&& (m+\al _{n})_{i,n-2}=1,
\qq(m+\beta _{qn})_{i,n-2}=1 \q  (q \geq i^{'}+1), \nonumber \\
&&(m+\beta _{qn})_{i,n-2}=0 \q  (q = i^{'}),
\qq(m+\beta _{qn}^{'})_{i,n-2}=1 \q  (q \geq i^{'}), \nonumber \\
&&(m+\gamma _{qn})_{i,n-2}=1 \q  (q \geq i^{'}+1),
\qq (m+\gamma _{qn}^{'})_{i,n-2}=1 \q  (q \geq i^{'}+1). \nonumber
\end{eqnarray}
From this, in particular, we also have the  
 Case 3) $1 \leq i \leq j =n-1$. \\ 
\q Case 4.1) $1 \leq i \leq j =n$ and $i=2i^{'}+1$:  \\
We use $e_{n-1}u=0$ and for any $m \in M_s$, 
\begin{eqnarray} 
&& (m+\al _{n-1})_{i,n}=1,    
\qq (m+\beta _{q,n-1})_{i,n}=1 \q  (q \geq i^{'}), \nonumber\\
&&(m+\beta _{q,n-1}^{'})_{i,n}=1 \q  (q \geq i^{'}+1),
\qq (m+\beta _{q,n-1}^{'})_{i,n}=0 \q  (q = i^{'}), \nonumber \\
&&(m+\gamma _{q,n-1})_{i,n}=1 \q  (q \geq i^{'}),
\qq (m+\gamma _{q,n-1}^{'})_{i,n}=1 \q  (q \geq i^{'}).\nonumber
\end{eqnarray}
\q Case 4.2) $1 \leq i \leq j =n$ and $i=2i^{'}$:  \\
We use $e_{n}u=0$ and for any $m \in M_s$, 
\begin{eqnarray} 
&& (m+\al _{n})_{i,n}=1,
\qq (m+\beta _{qn})_{i,n}=1 \q  (q \geq i^{'}), \nonumber \\
&& (m+\beta _{qn}^{'})_{i,n}=1 \q  (q \geq i^{'}+1),
\qq (m+\beta _{qn}^{'})_{i,n}=0 \q  (q = i^{'}), \nonumber \\
&& (m+\gamma _{qn})_{i,n}=1 \q  (q \geq i^{'}),
\qq (m+\gamma _{qn}^{'})_{i,n}=1 \q  (q \geq i^{'}). \nonumber
\end{eqnarray}
\q Case 5) $1 \leq j<i-1 < n-1$:  \\
We use $e_{i-1}u=0$ and for any $m \in M_s$, 
\begin{eqnarray} 
&& (m+\al _{i-1})_{ij}=1, 
\qq (m+\beta _{q,i-1})_{ij}=1 \q  (q \geq j), \no \\
&& (m+\beta _{q,i-1}^{'})_{ij}=1 \q  (q \geq j+1),
\qq (m+\beta _{q,i-1}^{'})_{ij}=0 \q  (q = j). \no 
\end{eqnarray}
\q Case 6) $1 \leq j=i-1 < n-1$:  
We use $e_{i}u=0$. \qed

\section{Irreducible $\uf$-module $\ue u_{\l ^{'}}(0)$} 
We keep the settings and notations as in Sect.2 and 3.
\subsection{Restricted specializations} 
In this subsection, 
we introduce the restricted specializations 
and its properties. 
\begin{df}
\label{Definition 3.1A}
Let $A:=\bbC [q,q^{-1}]$ and 
 $U_A^{res}$ be the $A$-subalgebra of $\uq$ generated by 
$\{ e_i^{(k)} , f_i^{(k)} , t_i^{\pm 1} | 
 i \in I , k \in \z_{+} \}$.
Let $l$ be an odd integer greater than $3$ and 
$\e$ be a primitive $l$-th root of unity such that 
$ \e ^{2d_i} \neq 1$ for any $i \in I$.
We regard $\bbC$ as $A$-algebra by
 $f(q).c:=f(\e) \cdot c$  for any $f(q) \in A, c \in \bbC$ 
and we denote it by $\bbC _{\e}$.
We define 
\[\ur := U_A^{res} \otimes _A \bbC _{\e},
\]
which is called  ``restricted specialization of $\ue$''.
Similarly, we define 
 $(\ur)^{+} , (\ur)^{-} , (\ur)^{0}$.
We denote $u \otimes 1$ as $u$ for any $u \in U_A^{res}$.
Let $\uf$ be the subalgebra of $\ur$ generated by 
$\{e_i , f_i , t_i^{\pm 1}\}_{i=1}^n$ 
(Similarly, 
  we define $(\uf)^{+}, (\uf)^{-}$ and $(\uf)^{0}$ ). 
\end{df}
 Next we review the representation theory of $\ur$. 
\begin{df}
\label{Definition 3.1B}
Let $L$ be a finite dimensional $\ur$-module.
If $t_i^l v=v$ for any $v \in V , i \in I$
(that is, $t_i^l$ is identity map),
we call $L$ `` $\ur$-module of type $1$''. 
\end{df}
\textbf{Remark.} (\cite{CP})
In general finite dimensional irreducible $\ur$-modules
are divided into $2^n$ types according to 
\{$\sigma : Q \longrightarrow \{\pm 1\}$ ; homomorphism of group \}).
Without a loss of generality, we may assume that  
finite dimensional irreducible $\ur$-modules are of type $1$.

\begin{df}
\label{Definition 3.1C}
For $\l =(\l_1 , \cdots , \l_n) \in \z _{+}^n$, 
let $I_{\l}$ be the left ideal of $\uq$ generated by 
$\{e_i , f_i^{\l_i +1} , t_i -q_i^{\l_i} | i \in I\}$
and $L(\l):=\uq / I_{\l}$.
We set $v_{\l}=1+I_{\l} \in L(\l)$.
Let $V_{A}^{res}(\l)$ be 
the $U_{A}^{res}$-submodule of $L(\l)$
generated by $v_{\l}$,
$V_{\e}^{res}(\l):=V_{A}^{res}(\l) \otimes _{A} \bbC_{\e}$,
and $W_{\e}^{res}(\l)$ be 
the maximal proper $\ur$-submodule of
$V_{\e}^{res}(\l)$.
We define
\[ L_{\e}^{res}(\l):=V_{\e}^{res}(\l) / W_{\e}^{res}(\l).
\]
\end{df}
\begin{thm}[\cite{L},\cite{L2}]
\label{Theorem 3.1A} 
\begin{enumerate}
\item
For any $\l \in \z_{+}^n$, $L_{\e}^{res}(\l)$
 is a finite dimensional irreducible $\ur$-module of type $1$
(We call $\l$ ``highest weight of $L_{\e}^{res}(\l)$''). 
\item
Let $L$ be a finite dimensional irreducible $\ur$-module
of type $1$. Then, there exists a unique element $\l \in 
\z_{+}^n$ such that $L \cong L_{\e}^{res}(\l)$. 
\item
Let $\l^{'}=(\l_1^{'}, \cdots , \l_n^{'}) \in \z_{l}^n 
\, (\z_l:=\{0,1, \cdots , l-1\})$,
$\l^{''} \in \z_{+}^n$, and $\l := \l^{'}+l\l^{''}$.
Then we have 
\[L_{\e}^{res}(\l) \cong 
 L_{\e}^{res}(\l^{'}) \otimes L_{\e}^{res}(l\l^{''}).
\]  
\end{enumerate}
\end{thm}
Next, we give the relation between the representations of 
$\ur$ and $\uf$. 
\begin{pro}[\cite{L}, \cite{L2}]
\label{Proposition 3.1A} 
\begin{enumerate}
\item
For any $\l =(\l_1 , \cdots , \l_n) \in \z_{l}^n$, 
we regard $L_{\e}^{res}(\l)$ as $\uf$-module and denote it 
$L_{\e}^{fin}(\l)$. 
Then $L_{\e}^{fin}(\l)$ is a 
finite dimensional irreducible $\uf$-module of type $1$.
Conversely, let $L$ be any finite dimensional irreducible
$\uf$-module of type $1$, 
then there exists a unique element 
 $\l =(\l_1 , \cdots , \l_n) \in \z_{l}^n$
such that $L \cong L_{\e}^{fin}(\l)$. 
\item
Let $U(\g)$ be the universal enveloping algebra of $\g$.
Then for any $\l \in \z_{+}^n$,
we can regard $L_{\e}^{res}(l\l)$ as 
a finite dimensional irreducible $U(\g)$-module of the highest weight $\l$. 
\end{enumerate}
\end{pro}

\subsection{Finite dimensional quantum algebra $\uf$}
In this subsection, we introduce the properties of $\uf$
to prove 
Theorem \ref{Theorem 3.3A}--\ref{Theorem 3.3C} below.
First, we introduce PBW theorem 
and the triangular decomposition of $\uf$.
\begin{thm}[\cite{L2}]
\label{Theorem 3.2A} 
Let $\beta_1 , \cdots , \beta _N$ be as in
Definition \ref{Definition 1.3A} 
($\Delta _{+}=\{\beta_1 , \cdots , \beta _N\}$).
We assume that $\g$ is not of type $G_2$.
Then we have 
\begin{enumerate}
\item
\{$e_{\beta _1}^{m_1} \cdots e_{\beta_N}^{m_N} |
 0 \leq m_i \leq l-1$ for any $1 \leq i \leq N$ \}
is a $\bbC$-basis of $(\uf)^{+}$. \\
\item
\{$f_{\beta _1}^{m_1} \cdots f_{\beta_N}^{m_N} |
 0 \leq m_i \leq l-1$ for any $1 \leq i \leq N$ \}
is a $\bbC$-basis of $(\uf)^{-}$. \\
\item
\{$t_{1}^{m_1} \cdots t_{n}^{m_n} |
 0 \leq m_i \leq 2l-1$ for any $1 \leq i \leq n$ \}
is a $\bbC$-basis of $(\uf)^{0}$ . 
\item
Let $\phi$ be the multiplication map 
\begin{eqnarray*}
\phi &:& (\uf)^{-} \otimes (\uf) ^{0} \otimes (\uf)^{+}
\longrightarrow \uf \\
&&\qq u_{-} \otimes u_{0} \otimes u_{+} \qq \mapsto 
\qq u_{-} u_0 u_{+}. 
\end{eqnarray*}
Then $\phi$ is an isomorphism of $\bbC$-vector space. 
\end{enumerate}
\end{thm}
By Theorem \ref{Theorem 3.2A}, we know that the dimension of $\uf$
is $2^n l^{n+2N}$. 
\begin{pro}[\cite{CP}]
\label{Proposition 3.2} 
We have $e_{\al}^l =f_{\al}^l =0$ in $\uf$ for any $\al \in 
\Delta _{+}$, and 
$t_i^{2l}=1$ in $\uf$ for any $i \in I$. 
\end{pro}
\begin{lem}
\label{Lemma 3.2}
We assume that $\g$ is not of type $G_2$. 
Let $\texttt{J}$ be the two-sided ideal of $\ue$ generated by 
$\{e_{\al}^l , f_{\al}^l | \al \in \Delta\}
\bigcup \{t_i^{2l}-1 | i \in I \}$. 
Then we have $\uf \cong \ue / \texttt{J}$. 
\end{lem}
{\sl Proof}.
 By the definition of $\uf$, $(e_i , f_i , t_i^{\pm 1})$
satisfies the relations in $\ue$.
Therefore, there exists the following 
$\bbC$-algebra homomorphism $\pi$ 
\[ \pi : \ue \longrightarrow \uf \q
(e_i , f_i , t_i^{\pm 1}) 
\mapsto (e_i , f_i , t_i^{\pm 1}).\]
In particular, by Proposition \ref{Proposition 3.2}, 
$\texttt{J} \subset Ker \pi$.  
 Conversely, by Theorem \ref{Theorem 1.3A}, 
for any $u \in \ue$, there exists  
$c(m)=c(m_1, \cdots , m_{2N+n} ) \in \bbC$
 $(m=(m_1, \cdots , m_{2N+n}) \in \z _{+}^{2N+n})$,
such that 
\[ u=\sum _{m \in \z_{+}^{2N+n}}c(m) 
f_{\beta_1}^{m_1} \cdots f_{\beta_N}^{m_N}
t_{1}^{m_{N+1}} \cdots t_{n}^{m_{N+n}}
e_{\beta_1}^{m_{N+n+1}} \cdots e_{\beta_N}^{m_{2N+n}}\, \in\, \ue.\] 
By Proposition \ref{pro1.2}, we have
\begin{eqnarray}
&& u \equiv\sum _{m \in \texttt{M}} \sum _{k_1, \cdots , k_n \geq 0}
c(m_1, \cdots , m_N , m_{N+1}+2k_1l, \cdots, m_{N+n}+2k_nl, \no \\ 
&&\qq m_{N+n+1}, \cdots m_{2N+n}) 
f_{\beta_1}^{m_1} \cdots f_{\beta_N}^{m_N}
t_{1}^{m_{N+1}} \cdots t_{_n}^{m_{N+n}}
e_{\beta_1}^{m_{N+n+1}} \cdots e_{\beta_N}^{m_{2N+n}}\, 
\q  \textrm{mod} \, \texttt{J}. \no 
\end{eqnarray}
where
\begin{eqnarray}
&& \texttt{M}:=\{m=(m_1,\cdots,m_{2N+n})\in \z_{+}^{2N+n} |
 0 \leq m_i <2l \q  
(N+1 \leq i \leq N+n), \no \\
&&\qq \qq 0 \leq m_i <l \q 
(1 \leq i \leq N , N+n+1 \leq i \leq 2N+n)\} \no 
\end{eqnarray}
\q Since $\pi (\texttt{J})=0$  by Proposition \ref{Proposition 3.2}, 
if $u \in Ker \pi$ then
\begin{eqnarray}
&& 0=\pi(u)=\sum _{m \in \texttt{M}} \sum _{k_1, \cdots , k_n \geq 0}
c(m_1, \cdots , m_N , m_{N+1}+2k_1l, \cdots, m_{N+n}+2k_nl, \no \\ 
&&\qq \qq m_{N+n+1}, \cdots m_{2N+n}) 
f_{\beta_1}^{m_1} \cdots f_{\beta_N}^{m_N}
t_{1}^{m_{N+1}} \cdots t_{_n}^{m_{N+n}}
e_{\beta_1}^{m_{N+n+1}} \cdots e_{\beta_N}^{m_{2N+n}}. \no 
\end{eqnarray}
\q Thus, by Theorem \ref{Theorem 3.2A},
\[  \sum _{k_1, \cdots , k_n \geq 0}
c(m_1, \cdots , m_N , m_{N+1}+2k_1l, \cdots, m_{N+n}+2k_nl, 
m_{N+n+1}, \cdots m_{2N+n})=0 \]
for any $m \in \texttt{M}$. 
Hence $u \in \texttt{J}$. \qed

\subsection{$\uf$-module structure on the 
$\ue u_{\l^{'}}(0)$}
In this subsection, we construct $\uf$-module
by using the Schnizer modules. 
\begin{lem}
\label{Lemma 3.3}
For $\g =\texttt{sp}(2n,\bbC)$ (resp. $\texttt{so}(2n+1,\bbC)$, 
$\texttt{so}(2n,\bbC)$)\q ($n \geq 2$), and
 $\l =(\l_1,\cdots,\l_n) \in \z^n$, 
let $a^{(0)}=(\e^{a_{ij}}),b^{(0)}=(\e^{b_{ij}})$ 
be as in Lemma \ref{Lemma 2.2A}, 
$(\Phi_{\l,a,b},V)$ be the Schnizer modules as in 
Theorem\ref{Theorem 2.1A}, 
and $u(0)=u(0,\cdots,0) \in V$ be the unique primitive vector.
Then we have
\begin{enumerate}
\item
$e_{\al}u(0)=0$\q for any $\al \in \Delta_{+}$.
\item
$t_i^{l}u(0)=u(0)$\q for any $i \in I.$
\item
$f_{\al}^lu(0)=0$\q for any $\al \in \Delta_{+}$.
\end{enumerate}
\end{lem}
{\sl Proof}.
(i) By Proposition \ref{Proposition 2.2A}, 
$e_iu(0)=0$\q for any $i \in I$.
On the other hand, by Proposition \ref{Proposition 1.3B}, 
$e_{\al} \in \ue^{+} \cap (\ue)_{\al}$\q for any $\al \in \Delta _{+}$.
Therefore $e_{\al}u(0)=0$. \\
\q (ii) By the explicit form of the action of $t_i$ 
in Theorem \ref{Theorem 2.1A},
there exists $c_{i} \in \z$ such that
\[t_iu(0)=\e^{\l_i+c_i}u(0).\]
Since $\l_i \in \z$, $t_i^{l}u(0)=u(0)$. \\
\q We shall prove (iii) in the next section. \qed 

We call  $\ue$-representation such that 
$e_i^l=f_i^l=0$ ``nilpotent representation''. By Lemma \ref{Lemma 3.2} 
(and proof of Lemma \ref{Lemma 3.3}), we can regard  
nilpotent irreducible $\ue$-representation (of type 1) as 
irreducible $\uf$-representation (of type 1).  
\begin{thm}
\label{Theorem 3.3A}
Let $\g=\texttt{sp}(2n,\bbC)\, (n \geq 2)$. 
For any $\l =(\l_1 , \cdots , \l_n) \in \z^n_l$, 
we define $\l^{'}:=(\l_1^{'}, \cdots , \l_n^{'}) \in \z^n$ by
\[ \l _j^{'}:=-\l_j-2 \q (1 \leq j \leq n-1),
\qq \l _n^{'}:=-2(\l_n+2).\]     
Let $a^{(0)}=(\e^{a_{ij}}),b^{(0)}=(\e^{b_{ij}})$ 
be as in Lemma \ref{Lemma 2.2A}, and
$(\Phi _{\l^{'},a,b}, V)$ be the Schnizer module
as in Theorem \ref{Theorem 2.1A}. 
We set $u_{\l^{'}}(0)=u(0, \cdots , 0)$, and
let $\ue u_{\l^{'}}(0)$
 be the $\ue$-submodule of $V$ generated by 
$u_{\l^{'}}(0)$.
Then we have 
\begin{enumerate}
\item
$\ue u_{\l^{'}}(0)$ is a $\uf$-module.\
\item
$\ue u_{\l^{'}}(0)$ is isomorphic to $L_{\e}^{fin}(\l)$
as $\uf$-module.
That is, $\ue u_{\l^{'}}(0)$ is a finite dimensional irreducible 
$\uf$-module of type $1$ with highest weight $\l$. 
\end{enumerate}
\end{thm}
{\sl Proof.}  
(i) Since $\Phi _{\l^{'},a,b} :\ue \longrightarrow End(V)$
is a homomorphism of $\bbC$-algebra,
$(\Phi _{\l^{'},a,b} (e_i),$ \,
$ \Phi _{\l^{'},a,b} (f_i),\Phi _{\l^{'},a,b}(t_i^{\pm 1}))$ 
satisfy the relations in Definition 2.1 in  $End(V)$.
On the other hand,
by Lemma \ref{Lemma 3.3} and Proposition \ref{pro1.2},
\[ (\Phi _{\l^{'},a,b}(e_i^l)) |_{\ue u_{\l^{'}}(0)}
=(\Phi _{\l^{'},a,b}(f_i^l)) |_{\ue u_{\l^{'}}(0)} =0
\q (\al \in \Delta _{+}),\]
\[ (\Phi _{\l^{'},a,b}(t_{i}^{\pm 2l})) |_{\ue u_{\l^{'}}(0) }=1
 \q (i \in I).\]
Therefore, by Lemma \ref{Lemma 3.2}, 
there exists a canonical homomorphism from
$\uf$ to $End(\ue u_{\l^{'}}(0) )$.
So we can regard $\ue u_{\l^{'}}(0)$ 
as $\uf$-module. \\
\q (ii)``Finite dimensionality'' of $\ue u_{\l^{'}(0)}$ is obvious. 
$\ue u_{\l^{'}}(0)$ is ``type 1'' by 
Lemma \ref{Lemma 3.3}(ii) and 
Proposition \ref{pro1.2}. 
So we shall prove the irreducibility of $\ue u_{\l^{'}}(0)$ 
 and that the highest weight of $\ue u_{\l^{'}}(0)$ is $\l$. \\
\q Irreducibility: 
We can also regard $\uf$ as well-defined $Q$-graded algebra 
by the following way (cf. Definition \ref{Definition 1.3B}).
\[ (\uf)_d :=\{ u+\texttt{J} | u \in (\ue)_d\} \q (d \in Q),\]
where $\texttt{J}$ is the two-sided ideal 
as in Lemma \ref{Lemma 3.2}.
Hence
\[ e_{i_1} \cdots e_{i_r} \in 
(\uf)_{\al_{i_1} +\cdots + \al_{i_r}}
\q (i_1, \cdots , i_r \in I).\]
On the other hand, by Proposition \ref{Proposition 3.2}, 
if $(\uf)_d \neq 0$ then 
\[d \leq (l-1)\sum _{\beta \in \Delta} \beta \]
 for any $d \in Q$, 
where $d \geq d^{'} \Leftrightarrow d-d^{'} \in Q_{+}$. 
So, there exists $r_0 \in \z_{+}$ such that
 $e_{i_1} e_{i_2} \cdots e_{i_r}=0$
for any $r \geq r_0$ and
$i_1, i_2 \cdots, i_r \in I$.
Thus, for any nonzero $\uf$-submodule $L$ of $\ue u_{\l^{'}}(0)$,
 there exists a nonzero element $v \in L$ such that 
$e_i v=0$ for any $i \in I$.
Therefore, by the uniqueness of primitive vector of 
Proposition \ref{Proposition 2.2A},
 $u_{\l^{'}}(0) \in L$.
Hence $L \supseteq \uf u_{\l^{'}}(0)
= \ue u_{\l^{'}}(0) \supseteq L$.  \\
\q Highest weight:
 By the definition of $L_{\e}^{fin}(\l)$, 
there exists a unique nonzero element (up to scalar multiplication)
$v \in L_{\e}^{fin}(\l)$ such that,
\[ e_i v=0 \q t_i v=\e_i^{\l_i}v \q \textrm{for any} \q i \in I,\] 
where $\e_i=\e^{d_i}$(since $\g=\texttt{sp}(2n,\bbC), 
(d_1, \cdots d_{n-1}, d_n)=(1, \cdots 1,2)$). 
So we shall prove that 
$t_i u_{\l^{'}}(0)=\e^{\l_i} u_{\l^{'}}(0)$ for any $i \in I$.
By the explicit form of $t_i$ in
Theorem \ref{Theorem 2.1A}, for any $i(1 \leq i \leq n-2)$,
we have
\begin{eqnarray}
&& T_{ij}u_{\l^{'}}(0) \no \\
&& =(\prod_{k=1}^{j-1}z_{k,j-1}^{-1}z_{kj}^2
z_{k,j+1}^{-1}z_{j+2,k}^{-1}z_{j+1,k}^2z_{jk}^{-1})\times
(z_{jj}^2z_{j,j+1}^{-1}z_{j+2,j}^{-1}z_{j+1,j}^2z_{j+1,j+1}^{-1}
z_{j+2,j+1}^{-1} \e ^{\l _j^{'}}) u_{\l^{'}}(0) \no \\
&&= \e^{c_j+\l_j^{'}}u_{\l^{'}}(0), \no
\end{eqnarray}
where
\begin{eqnarray}
 c_j&=&
\sum _{k=1}^{j-1}(-b_{k,j-1}+2b_{kj}-b_{k,j+1}-b_{j+2,k}
+2b_{j+1,k}-b_{jk}) \no \\
&& \qq \q + (2b_{jj}-b_{j,j+1}-b_{j+2,j}
+2b_{j+1,j}-b_{j+1,j+1}-b_{j+2,j+1}) \no \\
&=&\sum _{k=1}^{j-1}\{(b_{kj}-b_{k,j-1})-(b_{k,j+1}-b_{kj})\}
+\sum _{k=1}^{j-1}\{(b_{j+1,k}-b_{j+2,k})-(b_{jk}-b_{j+1,k})\} \no \\
&& \qq \q +(2b_{jj}-b_{j,j+1}-b_{j+1,j+1})
+2(b_{j+1,j}-b_{j+2,j})+(b_{j+2,j}-b_{j+2,j+1}) \no \\ 
&=&\sum _{k=1}^{j-1}(1-1)+\sum _{k=1}^{j-1}(1-1)
+(2-2-1)+2+1=2, \no
\end{eqnarray}
(since $k \leq j-1$, $b_{k,j-1}=b_{kj}-1=b_{k,j+1}-2$, 
$b_{jk}=b_{j+1,k}+1=b_{j+2,k}+2$).
Hence,
\[ t_j u_{\l^{'}}(0)=T_{1j}^{-1}u_{\l^{'}}(0)
=\e^{-2-\l_j^{'}}u_{\l^{'}}(0)=\e_i^{\l_j}u_{\l^{'}}(0).\]
Similarly, we can prove the case of $j=1,n-1,n$. \qed  \\
\q By the similar manner to the proof of Theorem \ref{Theorem 3.3A},
we can also prove the following theorems.
(if $\g=\texttt{so}(2n+1,\bbC)$ then 
$(d_1, \cdots d_{n-1}, d_n)=(2, \cdots ,2,1)$, and
 if $\g=\texttt{so}(2n,\bbC)$ then 
$(d_1, \cdots d_{n-1}, d_n)=(1, \cdots ,1,1)$).
\begin{thm}
\label{Theorem 3.3B}
Let $\g=\texttt{so}(2n+1,\bbC)\, (n \geq 3)$. 
For any $\l =(\l_1 , \cdots , \l_n) \in \z^n_l$, 
we denote $\l^{'}:=(\l_1^{'}, \cdots , \l_n^{'}) \in \z^n$ by
\[ \l _j^{'}:=-2(\l_j-2) \q (1 \leq j \leq n-1),
\qq  \l _n^{'}:=-\l_n-2.\]     
Let $a^{(0)}=(\e^{a_{ij}}),b^{(0)}=(\e^{b_{ij}})$ 
be in Lemma \ref{Lemma 2.2B}, and
$(\Phi _{\l^{'},a,b}, V)$ be the maximal cyclic $\ue$-representation 
of Theorem \ref{Theorem 2.1B}. 
We set $u_{\l^{'}}(0)=u(0, \cdots , 0)$.
Let $\ue u_{\l^{'}}(0)$
 be the $\ue$-submodule of $V$ generated by 
$u_{\l^{'}}(0)$.
Then we have 
\begin{enumerate}
\item
$\ue u_{\l^{'}}(0)$ is a $\uf$-module. 
\item
$\ue u_{\l^{'}}(0)$ is isomorphic to $L_{\e}^{fin}(\l)$
as $\uf$-module.
That is, $\ue u_{\l^{'}}(0)$ is a finite dimensional irreducible 
$\uf$-module of type $1$ with highest weight $\l$. 
\end{enumerate}
\end{thm}
\begin{thm}
\label{Theorem 3.3C}
Let $\g=\texttt{so}(2n,\bbC)\,(n \geq 4).$ 
For any $\l =(\l_1 , \cdots , \l_n) \in \z^n_l$, 
we denote $\l^{'}:=(\l_1^{'}, \cdots , \l_n^{'}) \in \z^n$ by
\[ \l _j^{'}:=-\l_j-2 \q (1 \leq j \leq n).\]
Let $a^{(0)}=(\e^{a_{ij}}),b^{(0)}=(\e^{b_{ij}})$ 
be in Lemma \ref{Lemma 2.2C}, and
$(\Phi _{\l^{'},a,b}, V)$ be 
the maximal cyclic $\ue$-representation 
 Theorem \ref{Theorem 2.1C}. 
We set $u_{\l^{'}}(0)=u(0, \cdots , 0)$.
Let $\ue u_{\l^{'}}(0)$
 be the $\ue$-submodule of $V$ generated by 
$u_{\l^{'}}(0)$.
Then we have 
\begin{enumerate}
\item
$\ue u_{\l^{'}}(0)$ is a $\uf$-module. 
\item
$\ue u_{\l^{'}}(0)$ is isomorphic to $L_{\e}^{fin}(\l)$
as $\uf$-module.
That is, $\ue u_{\l^{'}}(0)$ is a finite dimensional irreducible 
$\uf$-module of type $1$ with highest weight $\l$. 
\end{enumerate}
\end{thm}
\textbf{Comment}:
We expect that we can treat infinitesimal Verma modules
for orthogonal and symplectic cases by the similar way to
\cite{KN}.
\section{Proof of Lemma \ref{Lemma 3.3}{\rm (iii)} }
\renewcommand{\thesection}{\arabic{section}}
\setcounter{equation}{0}
\renewcommand{\theequation}{\thesection.\arabic{equation}}
\subsection{Case of $\l=(l-1, \cdots, l-1)$}
In this subsection we shall show Lemma \ref{Lemma 3.3} (iii)
for the special case :$\l=(l-1, \cdots, l-1)$.

For any $\l \in \z_l^n$, we denote the $\ue(\ge)$-module corresponding
to the representation 
$(\Phi_{a^{(0)},b^{(0)},\l^{'}},V)$ by $V(\l)$
($a^{(0)}, b^{(0)}$ as in Lemma \ref{Lemma 2.2A},
\ref{Lemma 2.2B}, \ref{Lemma 2.2C} and 
$\l^{'}$ as in 
Theorem \ref{Theorem 3.3A},
\ref{Theorem 3.3B}, \ref{Theorem 3.3C}.).
We obtain the following lemma by the similar way to the proof of 
Lemma \ref{Lemma 2.2A}. 
\begin{lem}
\label{f-act}
For $\ge=\texttt{sp}(2n,\bbC)$, 
the actions of $f_i$ on $V(\l)$ are given by the following formula. 
For any $m=(m_{i,j})_{i,j=1}^n \in M$, 
\begin{eqnarray*}
f_ju(m)=\sum_{i-1}^jv_{i,j}^m
\end{eqnarray*}
where
\begin{eqnarray*}
v_{i,j}^m&=&[m_{i,j}-m_{i,j+1}-m_{j+2,i}+2m_{j+1,i}-m_{j,i}+\mu_{i+1,j}^m]
u(m+\e_{i,j}) \\
&&+[m_{j+1,i}-m_{j,i}+\mu_{i+1,j}^m]u(m+\e_{j+1,i})
\q (1 \leq i <j \leq n-2), \\
v_{i,n-1}^m&=&[m_{i,n-1}-2m_{i,n}+2m_{n,i}-m_{n-1,i}+\mu_{i+1,n-1}^m]
u(m+\e_{i,n-1}) \\
&&+[m_{n,i}-m_{n-1,i}+\mu_{i+1,n-1}^m]u(m+\e_{n,i})
\q (1 \leq i <n-1), \\
v_{i,n}^m&=&[m_{i,n}-m_{n,i}+\mu_{i+1,n}^m]_{\e^2}u(m+\e_{i,n})
\q (1 \leq i <n), 
\end{eqnarray*}
\begin{eqnarray*}
v_{i,i}^m&=&[m_{i,i}-m_{i,i+1}-m_{i+2,i}+2m_{i+1,i}-m_{i+1,i+1}
-m_{i+2,i+1}-\l_i]
u(m+\e_{i,i}) \\
&&+[m_{i+1,i}-m_{i+1,i+1}-m_{i+2,i+1}-\l_i]u(m+\e_{j+1,i})
\q (1 \leq i <n-1), \\
v_{n-1,n-1}^m&=&[m_{n-1,n-1}-2m_{n-1,n}+2m_{n,n-1}-2m_{n,n}-\l_{n-1}]
u(m+\e_{n-1,n-1}) \\
&&+[m_{n,n-1}-2m_{n,n}-\l_{n-1}]u(m+\e_{n,n-1}), \\
v_{n,n}^m&=&[m_{n,n}-\l_n]_{\e^2}u(m+\e_{n,n}), 
\end{eqnarray*}
and
\begin{eqnarray*}
\mu_{i,j}=\sum_{k=i}^j\nu_{k,j}^m 
\q (1 \leq i \leq j \leq n),
\end{eqnarray*}
where
\begin{eqnarray*}
\nu_{i,j}^m&=&-m_{i,j-1}+2m_{i,j}-m_{i,j+1}-m_{j+2,i}+2m_{j+1,i}-m_{j,i} 
\q (1 \leq i <j<n-1),\\
\nu_{i,n-1}^m&=&-m_{i,n-2}+2m_{i,n-1}-2m_{i,n}-2m_{n,i}-m_{n-1,i} 
\q (1 \leq i <n-1),\\
\nu_{i,n}^m&=&-m_{i,n-1}+2m_{i,n}-m_{n,i}\q (1 \leq i <n),\\
\nu_{i,i}^m&=&2m_{i,i}-m_{i,i+1}-m_{i+2,i}+2m_{i+1,i}-m_{i+1,i+1}
-m_{i+2,i+1}-\l_i \q (1 \leq i <n-1),\\
\nu_{n-1,n-1}^m&=&2m_{n-1,n-1}-2m_{n-1,n}-2m_{n,n-1}-2m_{n,n}-\l_{n-1},\\
\nu_{n,n}^m&=&2m_{n,n}-\l_{n}.
\end{eqnarray*}
\end{lem}
\noindent Note that we can easily obtain the similar results for 
$\ge=\texttt{so}(m,\bbC)$ ($m=2n, 2n+1$).

For $\l=(\l_1, \cdots , \l_n) \in \z_l^n$, in case $\ge=\texttt{sp} (2n,\bbC)$,  
we define $m^{\l}=(m_{i,j}^{\l})_{i,j=1}^n \in M$ by 
\begin{eqnarray*}
m_{i,i}^{\l}&=&\l_i \q (1 \leq i \leq n), \\
m_{i,j}^{\l}&=&\l_i+ \cdots + \l_j \q (1 \leq i <j \leq n), \\
m_{j,i}^{\l}&=&\l_i+ \cdots + \l_{j-1} +2\l_j+ \cdots + 2\l_n
\q (1 \leq i <j \leq n). 
\end{eqnarray*}
Obviously, 
\begin{eqnarray}
&&m_{i,j+1}^{\l}-m_{i,j}^{\l}=\l_{j+1} 
\q (1 \leq i \leq j <n),\nn \\
&&m_{n,i}^{\l}-m_{i,n}^{\l}=\l_{n} \q (1 \leq i <n),
\label{mlm}\\
&&m_{j,i}^{\l}-m_{j+1,i}^{\l}=\l_{j} \q (1 \leq i \leq j <n-1).\nn
\end{eqnarray}
Similarly, for $\l=(\l_1, \cdots , \l_n) \in \z_l^n$, 
in case $\ge=\texttt{so}(2n+1,\bbC)$, we define
$m^{\l}=(m_{i,j}^{\l})_{i,j=1}^n \in M$ by 
\begin{eqnarray*}
m_{i,i}^{\l}&:=&\l_i 
\q (1 \leq i \leq n), \\
m_{i,j}^{\l}&:=&\l_i+ \cdots +\l_j 
\q (1 \leq i < j \leq n-1), \\
m_{i,n}^{\l}&:=&2\l_i +\cdots +2\l_{n-1}+\l_n 
\q (1 \leq i \leq n-1), \\
m_{j,i}^{\l}&:=&\l_i+ \cdots +\l_{j-1}+2\l_j+ \cdots +2\l_{n-1}+\l_n
\q (1 \leq i<j \leq n-1), \\
m_{n,i}^{\l}&:=&\l_i+ \cdots +\l_n 
\q (1 \leq i \leq n-1),
\end{eqnarray*}
and in case $\ge=\texttt{so}(2n,\bbC)$, we define
$m^{\l}=(m_{i,j}^{\l})_{1\leq i\leq n-1,1\leq j\leq n}^n \in M$ by 
\begin{eqnarray*}
m_{i,i}^{\l}&:=&\l_i 
\q (1 \leq i \leq n-1), \\
m_{n-1,n}^{\l}&:=&\l_n, \\
m_{i,j}^{\l}&:=&\l_i+ \cdots +\l_j 
\q (1 \leq i < j \leq n-2), \\
m_{i,n-1}^{\l}&:=&\l_i +\cdots +\l_{n-2}+\l_n 
\q (1 \leq i \leq n-2), 
\end{eqnarray*}
\begin{eqnarray*}
m_{i,n}^{\l}&:=&\l_i +\cdots +\l_{n-2}+\l_{n-1} 
\q (1 \leq i \leq n-2), \\
m_{j,i}^{\l}&:=&\l_i+ \cdots +\l_{j-1}+2\l_j+ \cdots +2\l_{n-2}+\l_{n-1}+\l_n
\q (1 \leq i<j \leq n-2), \\
m_{n-1,i}^{\l}&:=&\l_i+ \cdots +\l_n 
\q (1 \leq i \leq n-2).
\end{eqnarray*}

\begin{lem}
\label{lwv}
\begin{enumerate}
\item For any $\l \in \z_l^n$ and $\al \in \Delta_{+}$, we have
$f_{\al}u(m^{\l})=0$ in $V(\l)$. 
\item Set $\lm_0=(l-1, \cdots, l-1) \in \z_l^n$. 
For any $\al \in \Delta_{+}$ and $v \in V(\lm_0)$, 
write
\[
f_{\al}v=\sum_{m \in M} c_{\al}(m)u(m) \q (c_{\al}(m) \in \bbC).
\]
Then $c_{\al}(0)=0$, i.e. the vector $u(0)$ never occurs
in $f_{\al}v$.
\end{enumerate}
\end{lem}
{\sl Proof.}
We only show for $\ge=\texttt{sp}(2n,\bbC)$. 
The other cases are show similarly.
By Proposition \ref{Proposition 1.3B}, 
$f_{\al} \in \ue^{-}  \cap (\ue)_{\al}$. 
So, it is enough to prove the case of $\al=\al_i (i \in I)$. \\
(i) For any $1 \leq i<j <n-1$, 
by Lemma \ref{f-act} and (\ref{mlm}),
\begin{eqnarray*}
\nu_{i,j}^{m_{\l}}&=&-m_{i,j-1}^{\l}+2m_{i,j}^{\l}
-m_{i,j+1}^{\l}-m_{j+2,i}^{\l}+2m_{j+1,i}^{\l}-m_{j,i}^{\l} \\
&=&(m_{i,j}^{\l}-m_{i,j-1}^{\l})-(m_{i,j+1}^{\l}-m_{i,j}^{\l})
+(m_{j+1,i}^{\l}-m_{j+2,i}^{\l})-(m_{j,i}^{\l}-m_{j+1,i}^{\l}) \\
&=&\l_j-\l_{j+1}+\l_{j+1}-\l_j=0. 
\end{eqnarray*}
Similarly, we obtain 
\begin{eqnarray*}
\nu_{i,n-1}^{m_{\l}}=0 \q  (1 \leq i <n-1), 
\qq \nu_{i,n}^{m_{\l}}=0  \q  (1 \leq i <n).
\end{eqnarray*}
Further, for any $1 \leq i <n-1$,
\begin{eqnarray*}
\nu_{i,i}^{m_{\l}}&=&2m_{i,i}^{\l}-m_{i,i+1}^{\l}-m_{i+1,i+1}^{\l}
-\l_i-m_{i+2,i}^{\l}+2m_{i+1,i}^{\l}-m_{i+2,i+1}^{\l}\\
&=& 2\l_i-(\l_i+\l_{i+1})-\l_{i+1}-\l_{i}
-(\l_i+\l_{i+1}+2\l_{i+2}+ \cdots 2\l_n)\\
&&\q +2(\l_i+2\l_{i+1}+2\l_{i+2}+ \cdots 2\l_n)
-(\l_{i+1}+2\l_{i+2}+ \cdots 2\l_n) \\
&=&\l_i
\end{eqnarray*}
Similarly, we obtain 
\begin{eqnarray*}
\nu_{n-1,n-1}^{m_{\l}}=\l_{n-1}, 
\q \nu_{n,n}^{m_{\l}}=\l_n.
\end{eqnarray*}
Thus, it follows from Lemma \ref{f-act} that
 for any $1 \leq i,j \leq n$, 
\begin{eqnarray*}
\mu_{i,j}^{m^{\l}}=\sum_{k=i}^j
 \nu_{k,j}^{m^{\l}}=\nu_{j,j}^{m^{\l}}=\l_j.
\end{eqnarray*}
Hence, for any $1 \leq i<j < n-1$,
\begin{eqnarray*}
v_{i,j}^{m^{\l}}&=&
[-(m_{i,j+1}^{\l}-m_{i,j}^{\l})+(m_{j+1,i}^{\l}-m_{j+2,i}^{\l})
-(m_{j,i}^{\l}-m_{j+1,i}^{\l})-\mu_{i+1,j}^{m^{\l}}]u(m^{\l}+\e_{i,j}) \\
&& \q+[-(m_{j,i}^{\l}-m_{j+1,i})-\mu_{i+1,j}^{m^{\l}}]u(m^{\l}+\e_{j+1,i}) \\
&=&[-\l_{j+1}+\l_{j+1}-\l_j+\l_j]u(m^{\l}+\e_{i,j})
+[-\l_j+\l_j]u(m^{\l}+\e_{j+1,i})=0.
\end{eqnarray*}
Similarly, we obtain 
\begin{eqnarray*}
 v_{i,j}^{m^{\l}}=0 \q (1 \leq i \leq j \leq n).
\end{eqnarray*}
Therefore, by Lemma \ref{f-act}, for any $j \in I$ we have 
$f_ju(m^{\l})=\sum_{i=1}^j v_{i,j}^{m^{\l}}=0.$\\
(ii) For any $\l \in \z_l^n$, 
$v=\sum_{m \in M} c(m)u(m) \in V(\l)$, by Lemma \ref{f-act}, 
\begin{eqnarray*}
f_jv=\sum_{m \in M} c(m)(\sum_{i=1}^j v_{i,j}^m)
=\sum_{i=1}^j(\sum_{m \in M}c(m) v_{i,j}^m) \q (j \in I).
\end{eqnarray*}
Since for any $1 \leq i <j<n-1$, $m_{i,j}$ does not appear in 
$\mu_{i+1,j}^m$ and $(i,j)\ne(i,j+1), (j+2,i), (j+1,i), (j,i)$,
we have 
\begin{eqnarray*}
\sum_{m \in M}c(m) v_{i,j}^m &=&\sum_{m \in M}c(m) 
[m_{i,j}-m_{i,j+1}-m_{j+2,i}+2m_{j+1,i}-m_{j,i}+\mu_{i+1,j}^m]
u(m+\e_{i,j})\\
&=&\sum_{m \in M}c(m-\e_{i,j}) 
[(m_{i,j}-1)-m_{i,j+1}-m_{j+2,i}+2m_{j+1,i}-m_{j,i}+\mu_{i+1,j}^m]u(m)
\label{cmv}
\end{eqnarray*}
If $m\equiv 0$, we have $\mu_{i+1,j}^m=-\lm_j$. Thus, 
we obtain that the coefficient of  $u(0)$ is equal to:
\begin{eqnarray*}
c(-\e_{i,j})[-1-\l_j].
\end{eqnarray*}
Hence, if $\l=(l-1, \cdots , l-1)$, then this is $0$. 
Similarly, we obtain that $u(0)$ does not
appear in 
$\sum _{m \in M}c(m) v_{i,j}^m$ for all other $i,j \in I$. 
So the coefficient of $u(0)$ in $f_jv$ is equal to $0$.
\qed
\begin{lem}
\label{-1-case}
We have 
$f_{\al}^lu(0)=0$ in $V(\lm_0)$ 
$(\lm_0=(l-1,\cdots,l-1))$ for any $\al \in \Delta_{+}$.
\end{lem}
{\sl Proof.} 
By Proposition \ref{pro1.2}, for any $\al \in \Delta_{+}$, 
$f_{\al}^l$ is a central element of $\ue$. 
Thus, 
\begin{eqnarray*}
e_i(f_{\al}^lu(0))=f_{\al}^l(e_iu(0))=0 
\q (i \in I).
\end{eqnarray*}
So, $f_{\al}^lu(0)$ is a primitive vector. 
Therefore, by the uniqueness of primitive vector 
(see Proposition \ref{Proposition 2.2A}), 
$f_{\al}^lu(0) \in \bbC u(0)$. \\
\q On the other hand, by Lemma \ref{lwv}(ii), 
the coefficient of $u(0)$ in  $f_{\al}^lu(0)$ is $0$. 
Hence $f_{\al}^lu(0)=0$. \qed 
\subsection{Proof of $e_{\al}^l=0$ on $V(\l)$}
\begin{df}
\label{VM}
Let $\l=(\l_1, \cdots, \l_n) \in \z_l^n$, 
$I_{\l}$ be the left ideal of $\ue$ generated by 
$\{e_i, t_i-\e^{\l_i}_i, f_{\al}^l\,
|\, i \in I, \al \in \Delta_{+}\}$.
We set $M(\l):=\ue / I_{\l}$.
\end{df}
\begin{pro}
\label{VM-pro}
(\cite{DK} Proposition 3.2, Corollary 3.2(b)) \\
\q (i) If $\l=(l-1, \cdots, l-1)$, then $M(\l)$ is 
an irreducible $\ue$-module. \\
\q (ii) For any $\l \in \z_l^n$, 
$\textrm{dim}M(\l)=l^{n^2} \q (=\textrm{dim}V(\l))$.
\end{pro}
\begin{pro}
\label{iso}
For $\l=(l-1, \cdots, l-1)$, 
$M(\l) \cong V(\l)$ (as $\ue$-module).
\end{pro}
{\sl Proof.} 
By Lemma \ref{-1-case} and the property of $u(0)$, we have
\begin{eqnarray*}
e_iu(0)=0, \q t_iu(0)=\e_{i}^{\l_i}u(0), 
\q f_{\al}^lu(0)=0  
\q (i \in I, \al \in \Delta_{+}). 
\end{eqnarray*}
So, by the universality of $M(\l)$, 
there exists an $\ue$-module homomorphism 
$\phi : M(\l) \longrightarrow V(\l)$ such that 
$\phi(1+I_{\l})=u(0)$. 
By Proposition \ref{VM-pro}(i), 
$M(\l)$ is an irreducible $\ue$-module if 
$\l=(l-1, \cdots, l-1)$, 
and $\phi\not\equiv 0$. 
Hence $\phi$ is injective. 

On the other hand, by Proposition \ref{VM-pro}(ii), 
$\textrm{dim}M(\l)=\textrm{dim}V(\l)$. 
Thus $\phi$ is surjective. 
Therefore $\phi$ is an isomorphism of $\ue$-module. 
\qed
\begin{lem}
\label{e=0}
For any $\l \in \z_l^n$ and $\al \in \Delta_{+}$, 
$e_{\al}^l=0$ on $V(\l)$.
\end{lem}
{\sl Proof.} 
By Proposition \ref{iso}, Proposition \ref{VM-pro}(i), 
$V(l-1, \cdots, l-1)$ is an irreducible $\ue$-module and then
we have $\ue u(0)=V(l-1, \cdots, l-1)$. 
Thus, 
\begin{eqnarray*}
e_{\al}^lV(l-1, \cdots, l-1)=e_{\al}^l(\ue u(0))
=\ue(e_{\al}^lu(0))=\{0\}.
\end{eqnarray*}
Hence $e_{\al}^l=0$ on $V(l-1, \cdots, l-1)$. 
Due to Lemma \ref{Lemma 2.2A} we know that
the actions of $e_{i}$ on $V(\l)$ do not depend on $\l$. 
Therefore, for any $\l \in \z_l^n$, 
$e_{\al}^l=0$ on $V(\l)$. 
\qed

\subsection{General case}
\begin{lem}
\label{uv=u0}
For any $\l \in \z_l^n$ and $v \in V(\l) (v \neq 0)$, 
there exists $u^{+} \in \ue^{+}$ 
such that $u^{+}v=u(0)$.
\end{lem}
{\sl Proof.} 
By Lemma \ref{e=0}, we can regard $V(\l)$ as a
$(\uf)^{+}$-module. 
(see proof of Lemma \ref{Lemma 3.2}). 
So, by the similar manner to the proof of Theorem \ref{Theorem
3.3A}(ii), we can take $u^{+} \in \ue^{+}$ such that 
$u^{+}v$ is a nonzero primitive vector. 
Therefore, by the uniqueness of the primitive vector, 
we have $u^{+}v \in \bbC^{\times} u(0)$. 
\qed

{\sl Proof of Lemma \ref{Lemma 3.3}(iii)}. Let us show that 
for any $\l \in \z_l^n$ and $\al \in \Delta_{+}$, 
$f_{\al}^lu(0)=0$ in $V(\l)$.
By Lemma \ref{uv=u0}, 
there exists $u^{+} \in \ue^{+}$ such that 
$u^{+}u(m^{\l})=u(0)$. 
Since $f_{\al}u(m^{\l})=0$ by Lemma \ref{lwv}(i), 
\begin{eqnarray*}
f_{\al}^lu(0)=f_{\al}(u^{+}u(m^{\l}))=u^{+}(f_{\al}^lu(m^{\l}))=0
\end{eqnarray*}
\qed


\end{document}